\documentclass[12pt]{amsart}
\usepackage{amsmath}
\usepackage{amsfonts}
\begin{document}
\numberwithin{equation}{section}

\def\1#1{\overline{#1}}
\def\2#1{\widetilde{#1}}
\def\3#1{\widehat{#1}}
\def\4#1{\mathbb{#1}}
\def\5#1{\frak{#1}}
\def\6#1{{\mathcal{#1}}}

\newcommand{\de}{\partial}
\newcommand{\R}{\mathbb R}
\newcommand{\al}{\alpha}
\newcommand{\tr}{\widetilde{\rho}}
\newcommand{\tz}{\widetilde{\zeta}}
\newcommand{\tv}{\widetilde{\varphi}}
\newcommand{\hv}{\hat{\varphi}}
\newcommand{\tu}{\tilde{u}}
\newcommand{\tF}{\tilde{F}}
\newcommand{\debar}{\overline{\de}}
\newcommand{\Z}{\mathbb Z}
\newcommand{\C}{\mathbb C}
\newcommand{\Po}{\mathbb P}
\newcommand{\zbar}{\overline{z}}
\newcommand{\G}{\mathcal{G}}
\newcommand{\So}{\mathcal{S}}
\newcommand{\Ko}{\mathcal{K}}
\newcommand{\U}{\mathcal{U}}
\newcommand{\B}{\mathbb B}
\newcommand{\oB}{\overline{\mathbb B}}
\newcommand{\Cur}{\mathcal D}
\newcommand{\Dis}{\mathcal Dis}
\newcommand{\Levi}{\mathcal L}
\newcommand{\SP}{\mathcal SP}
\newcommand{\Sp}{\mathcal Q}
\newcommand{\A}{\mathcal O^{k+\alpha}(\overline{\mathbb D},\C^n)}
\newcommand{\CA}{\mathcal C^{k+\alpha}(\de{\mathbb D},\C^n)}
\newcommand{\Ma}{\mathcal M}
\newcommand{\Ac}{\mathcal O^{k+\alpha}(\overline{\mathbb D},\C^{n}\times\C^{n-1})}
\newcommand{\Acc}{\mathcal O^{k-1+\alpha}(\overline{\mathbb D},\C)}
\newcommand{\Acr}{\mathcal O^{k+\alpha}(\overline{\mathbb D},\R^{n})}
\newcommand{\Co}{\mathcal C}
\newcommand{\Hol}{{\sf Hol}(\mathbb H, \mathbb C)}
\newcommand{\Aut}{{\sf Aut}(\mathbb D)}
\newcommand{\D}{\mathbb D}
\newcommand{\oD}{\overline{\mathbb D}}
\newcommand{\oX}{\overline{X}}
\newcommand{\loc}{L^1_{\rm{loc}}}
\newcommand{\la}{\langle}
\newcommand{\ra}{\rangle}
\newcommand{\thh}{\tilde{h}}
\newcommand{\N}{\mathbb N}
\newcommand{\kd}{\kappa_D}
\newcommand{\Hr}{\mathbb H}
\newcommand{\ps}{{\sf Psh}}
\newcommand{\Hess}{{\sf Hess}}
\newcommand{\subh}{{\sf subh}}
\newcommand{\harm}{{\sf harm}}
\newcommand{\ph}{{\sf Ph}}
\newcommand{\tl}{\tilde{\lambda}}
\newcommand{\gdot}{\stackrel{\cdot}{g}}
\newcommand{\gddot}{\stackrel{\cdot\cdot}{g}}
\newcommand{\fdot}{\stackrel{\cdot}{f}}
\newcommand{\fddot}{\stackrel{\cdot\cdot}{f}}
\def\v{\varphi}
\def\Re{{\sf Re}\,}
\def\Im{{\sf Im}\,}

\title[Pluripotential theory and infinitesimal generators]{Pluripotential theory,
semigroups and boundary behavior of infinitesimal generators in
strongly convex domains}
\author[F. Bracci]{Filippo Bracci}
\address{F. Bracci: Dipartimento Di Matematica\\
Universit\`{a} di Roma \textquotedblleft Tor Vergata\textquotedblright\ \\
Via Della Ricerca Scientifica 1, 00133 \\
Roma, Italy} \email{fbracci@mat.uniroma2.it}
\thanks{}
\author[M.D. Contreras]{Manuel D. Contreras$^\dag$ }
\address{M.D. Contreras \and  S. D\'{\i}az-Madrigal: Camino de los Descubrimientos, s/n\\
Departamento de Matem\'{a}tica Aplicada II \\
Escuela T\'ecnica Superior de Ingenieros\\
Universidad de Sevilla\\
41092, Sevilla\\
Spain.} \email{contreras@esi.us.es, madrigal@us.es}
\author[S. D\'{\i}az-Madrigal]{Santiago D\'{\i}az-Madrigal$^\dag$}

\date{26 July 2006}



\thanks{$^\dag$Partially supported by the \textit{Ministerio
de Ciencia y Tecnolog\'{\i}a} and the European Union (FEDER) project
BFM2003-07294-C02-02 and by \textit{La Consejer\'{\i}a de
Educaci\'{o}n y Ciencia de la Junta de Andaluc\'{\i}a.}}

\begin{abstract} We characterize infinitesimal generators of
semigroups of holomorphic self-maps of strongly convex domains
using the pluricomplex Green function and the pluricomplex Poisson
kernel. Moreover, we study boundary regular fixed points of
semigroups. Among other things, we characterize boundary regular
fixed points both in terms of the boundary behavior of
infinitesimal generators and in terms of pluripotential theory.
\end{abstract}

\maketitle

\tableofcontents

\def\Label#1{\label{#1}}


\def\cn{{\C^n}}
\def\cnn{{\C^{n'}}}
\def\ocn{\2{\C^n}}
\def\ocnn{\2{\C^{n'}}}
\def\je{{\6J}}
\def\jep{{\6J}_{p,p'}}
\def\th{\tilde{h}}


\def\dist{{\rm dist}}
\def\const{{\rm const}}
\def\rk{{\rm rank\,}}
\def\id{{\sf id}}
\def\aut{{\sf aut}}
\def\Aut{{\sf Aut}}
\def\CR{{\rm CR}}
\def\GL{{\sf GL}}
\def\Re{{\sf Re}\,}
\def\Im{{\sf Im}\,}
\def\U{{\sf U}}

\def\la{\langle}
\def\ra{\rangle}

\emergencystretch15pt \frenchspacing

\newtheorem{Thm}{Theorem}[section]
\newtheorem{Cor}[Thm]{Corollary}
\newtheorem{Pro}[Thm]{Proposition}
\newtheorem{Lem}[Thm]{Lemma}

\theoremstyle{definition}\newtheorem{Def}[Thm]{Definition}

\theoremstyle{remark}
\newtheorem{Rem}[Thm]{Remark}
\newtheorem{Exa}[Thm]{Example}
\newtheorem{Exs}[Thm]{Examples}

\def\bl{\begin{Lem}}
\def\el{\end{Lem}}
\def\bp{\begin{Pro}}
\def\ep{\end{Pro}}
\def\bt{\begin{Thm}}
\def\et{\end{Thm}}
\def\bc{\begin{Cor}}
\def\ec{\end{Cor}}
\def\bd{\begin{Def}}
\def\ed{\end{Def}}
\def\br{\begin{Rem}}
\def\er{\end{Rem}}
\def\be{\begin{Exa}}
\def\ee{\end{Exa}}
\def\bpf{\begin{proof}}
\def\epf{\end{proof}}
\def\ben{\begin{enumerate}}
\def\een{\end{enumerate}}

\section*{Introduction}

A (continuous) semigroup $(\Phi_t)$ of holomorphic functions in a
domain $D\subset\C^n$ is a continuous homomorphism from the additive
semigroup of non-negative real numbers into the composition
semigroup of all holomorphic self-maps of $D$  endowed with the
compact-open topology. Namely, the map $[0,+\infty )\ni t\mapsto
(\Phi_{t})\in\mathrm{Hol}(D,D)$ satisfies the following conditions:

\begin{enumerate}
\item $\Phi _{0}$ is the identity map $\id_D$ in $D,$

\item $\Phi _{t+s}=\Phi _{t}\circ \Phi _{s},$ for all
$t,s\geq 0,$

\item $\Phi _{t}$ tends to $\id_D$ as $t$ tends to $0$
uniformly on compacta of $D$.
\end{enumerate}

It is well known after the basic work of Berkson and Porta
\cite{Berkson-Porta} in the unit disc $\D=\{\zeta\in\C: |\zeta|<1\}$
 that the dependence of   every
 semigroup $(\Phi_t)$ of holomorphic self-maps of a domain
$D\subset\C^n$ on the variable $t$ is analytic and to each
continuous semigroup $(\Phi_t)$ there corresponds a holomorphic
vector field $F:D\to \C^n$ such that $\frac{\de \Phi_t}{\de
t}=F\circ \Phi_t$. This vector field $F$ is called the {\sl
infinitesimal generator} of the semigroup $(\Phi_t)$. Conversely,
if a holomorphic vector field $F:D\to \C^n$ is {\sl semicomplete},
namely for all $z\in D$ its local flow $\gamma^z(t)$ such that
$\gamma^z(0)=z$ is well defined for all $t\geq 0$, then $F$ is the
infinitesimal generator of a
 semigroup of holomorphic self-maps of $D$. We refer to
\cite[Section 2.5.3]{Abate}  and \cite{R-S2} for more details. Be
aware that in the literature there is not a standard sign convention
for the Cauchy problem generating $F$, namely, sometimes the problem
$\frac{\de \Phi_t}{\de t}=-F\circ \Phi_t$ is considered and thus all
formulas regarding $F$ have reverse inequalities with respect to our
formulas. For instance, regarding the bibliography of the present
paper, such a convention is adopted in \cite{AERS}, \cite{A-S},
\cite{E-R-S}, \cite{E-S}, \cite{E-S2}, \cite{R-S}, \cite{R-S2} and
\cite{Shoiket}.

It is clear that the analytical properties of an infinitesimal
generator are strictly related to the dynamical and geometrical
properties of its semigroup. For instance, any zero of $F$ in $D$
corresponds to a common fixed point for $(\Phi_t)$.

Therefore one of the main questions in the theory of semigroups of
holomorphic functions is that of characterizing (in the most
useful way) those holomorphic vector fields which are
infinitesimal generators. For $D=\D$, the unit disc of $\C$, there
is a very nice representation formula, due to Berkson and Porta
\cite{Berkson-Porta} (see also \cite{Abate} and \cite{Shoiket}).
Namely:

\bt[Berkson-Porta]\label{24} A  holomorphic function
$G:\mathbb{D\rightarrow C}$ is the infinitesimal generator of a
semigroup $(\Phi_t)$ in $\mathbb{D}$ if and only if there exists a
point $b\in \overline{\mathbb{D}}$ and a holomorphic function
$p:\mathbb{D\rightarrow C}$ with $\Re p\geq 0$ such that
\begin{equation*}
G(z)=(z-b)(\overline{b}z-1)p(z),\text{ \quad }z\in \mathbb{D}.
\end{equation*}
\et

If the semigroup is not an elliptic group (that is, some/all
iterates $\Phi_t$  for $t>0$  are elliptic automorphisms), the point
$b$ given in Berkson-Porta's formula is exactly the {\sl
Denjoy-Wolff point} of the semigroup $(\Phi_t)$. Namely, $\lim_{t\to
\infty} \Phi_t(z)=b$ for all $z\in \D$ (see also Section two). Other
alternative descriptions of infinitesimal generators in $\D$ can be
found in \cite[Section 3.6]{Shoiket}.

In several variables there are various characterizations of
infinitesimal generators (see \cite{R-S2} for a good account). All
these characterizations reflect the basic fact that holomorphic
self-maps of a domain are contractions for the Kobayashi metric
 of such a domain. In fact, Abate \cite{Abate2}
proved that if $D$ is a strongly convex domain with smooth boundary
and with Kobayashi metric $k_D$, then a holomorphic vector field
$F:D\to\C^n$ is an infinitesimal generator if and only if
$d(k_D\circ F)(z)\cdot F(z)\leq 0$ for all $z\in D$. Unfortunately,
even for the case $D=\B^n$, the unit ball of $\C^n$, such a formula
is rather complicated and does not give any information on the
dynamical properties of the associated semigroup. Later on, still in
this optic, C. de Fabritiis  gave a better description of a class of
infinitesimal generators called ``of one-dimensional type'' (see
\cite{deF}). Some rather precise characterizations of infinitesimal
generators in the unit ball of infinite dimensional Hilbert spaces
are given by D. Aharonov, M. Elin, S. Reich, and D. Shoikhet in
\cite{AERS}, \cite{E-S2} and \cite{R-S}.

Part of the present paper is devoted to find characterizations of
infinitesimal generators in bounded strongly convex domains with
smooth boundary (here and in the rest of the paper ``smooth''
means at least of class $C^3$) by means of the pluricomplex Green
function $G_D$ of Klimek \cite{Kl}, Lempert \cite{Le} and Demailly
\cite{De} and the pluricomplex Poisson kernel $u_{D,p}$ introduced
by Patrizio and the first named author in \cite{B-P} (see Section
one for definitions and preliminaries about pluripotential theory
in strongly convex domains). In particular, we prove (see Theorems
\ref{caratterizza-Green} and \ref{charnobene}):

\bt Let $F:D\to \C^n$ be a holomorphic vector field. The following
are equivalent:
\begin{enumerate}
\item  The map $F$ is an infinitesimal generator  of a semigroup
of holomorphic self-maps of $D$.
\item For all $z,w\in D$, $z\neq
w$, it holds $d(k_D)|_{(z,w)}\cdot (F(z), F(w))\leq 0$.
\item  For
all $z,w\in D$, $z\neq w$, it holds $d(G_D)_{(z,w)}\cdot
(F(z),F(w))\leq 0$.
\item For all  $z,w\in D$ and for all $r>0$
such that $z-rF(z),w-rF(w)\in D$ it holds
$k_{D}(z-rF(z),w-rF(w))\geq k_{D}(z,w)$.
\end{enumerate}
Moreover, if $F$ is $C^1$-regular at a point $p\in\de D$, then $F$
is an infinitesimal generator whose associated semigroup has
Denjoy-Wolff point at $p\in\de D$ if and only if $d(u_{D,p})_z\cdot
F(z)\leq 0$ for all $z\in D$. \et

In case $D=\B^n$ is the unit ball of $\C^n$ (or more generally for
the unit ball of complex Hilbert spaces), equivalence between (1)
and (4)  (and also with an explicit expression of (3), see Remark
\ref{sciocco}) was proven with different methods by Reich and
Shoikhet \cite[Theorem 2.1]{R-S}. The last statement can be seen as
a Berkson-Porta like formula at the boundary. Moreover, this last
formula is just a particular case of a general one for the existence
of {\sl boundary regular fixed points}. We recall that a point
$p\in\de D$ is a boundary regular fixed point---BRFP for short---for
a semigroup $(\Phi_t)$ if it is a fixed point for non-tangential
limits for all $\Phi_t$'s and if the {\sl boundary dilatation
coefficients} at $p$ of the $\Phi_t$'s are all finite (roughly
speaking, the boundary dilatation coefficient of a self-map $f$ of
$D$ at $p$ is a measure of the velocity $f$ approaches $p$ when
moving to $p$; see Section two for details and precise definitions).

The second part of this paper is devoted to characterize BRFPs of
semigroups in terms of the pluricomplex Poisson kernel $u_{D,p}$
and the local behavior of the infinitesimal generator. In this
direction, we cite the following result from \cite{E-S2} (see also
\cite{E-R-S}, \cite{A-S})

\bt[Elin-Shoikhet]\label{Shoikhet} Let $F:\B^n\to \C^n$ be the
infinitesimal generator of a semigroup $(\Phi_t)$ in $\B^n$ and
$p\in\de\B^n$. Assume that $\lim_{(0,1)\ni r\to 1}F(r p)=0$. The
following are equivalent:
\begin{enumerate}
\item $\liminf_{(0,1)\ni r\to 1}\Re \la F(r
p),p\ra/(r-1)<+\infty$. \item $\lim_{(0,1)\ni r\to 1}\la F(r
p),p\ra/(r-1)=\beta$ exists finitely. \item The point $p$ is a
BRFP for the semigroup $(\Phi_t)$.
\end{enumerate}
Moreover, if one of the three conditions holds,  then $\beta\in \R
$ and the boundary dilatation coefficient  of $\Phi_t$ at $p$ is
$e^{t\beta}$. \et

The hypothesis in Theorem \ref{Shoikhet} that the infinitesimal
generator $F$ has radial limit $0$ at $p$, being essential in the
proof of their result, is however not necessary for a point $p$ to
be a BRFP (see Example \ref{unos}). Moreover, and surprisingly
enough, Theorem \ref{Shoikhet} would be false without such an
hypothesis (see Example \ref{dues} where it is constructed an
infinitesimal generator for which (1) holds at some $p\in\de \B^2$
but $p$ is not a BRFP for the associated semigroup). In fact, it
turns out that a point $p\in \de \B^n$ is a BRFP for the semigroup
if and only if a condition  similar to (1) holds not just for the
radial direction but for all the directions. To be more precise and
in order to state the result for general strongly convex domains, we
need to use the so called {\sl Lempert projection devices}. For the
time being, we can say that a Lempert projection device $(\v,
\tr_\v)$ is given by a particular holomorphic map $\v:\D\to D$
(called {\sl complex geodesic}) which extends smoothly on $\de \D$
and a holomorphic map $\tr_\v:D\to \D$ such that $\tr_\v\circ
\v={\sf id}_\D$ (actually a Lempert projection device is a triple of
maps,  we refer the reader to Section one for details). For the unit
ball $\B^n$ a Lempert projection device $(\v, \tr_\v)$ is nothing
but a (suitable) parametrization $\v:\D\to \B^n$ of the intersection
of $\B^n$ with an affine complex line and $\tr_\v$ is the orthogonal
projection  on it (see also Section six where the case of $\B^n$ is
studied in detail). Our second main result is the following:

\bt Let $D\subset \C^n$ be a bounded strongly convex domain with
smooth boundary, let $F$ be the infinitesimal generator of a
semigroup $(\Phi_t)$ of holomorphic self-maps of $D$ and $p\in\de
D$. The following are equivalent:
\begin{enumerate}
\item The
semigroup $(\Phi_t)$ has a BRFP at $p$ with boundary dilatation
coefficients $\al_t(p)\leq e^{\beta t}$ for all $t\geq 0$.
\item There exists $\beta\in \R$ such that $d( u_{D,p})_z\cdot F(z)
+\beta u_{D,p}(z)\leq 0$ for all $z\in D$.
\item There exists $C>0$ such that for any Lempert's projection
device $(\v, \tr_v)$ with $\v(1)=p$  it follows
\[
\limsup_{(0,1)\ni r\to 1}\frac{|d(\tr_\v)_{\v(r)} \cdot
F(\v(r))|}{1-r}\leq C.
\]
\end{enumerate}
Moreover, if $p$ is a BRFP for $(\Phi_t)$ with boundary dilatation
coefficients $\al_t(p)=e^{-b t}$ then
\[
b=\inf_{z\in D} \frac{d(u_{D,p})_z\cdot F(z)}{u_{d,p}(z)},
\]
and  the non-tangential limit
\[
A(\v,p):=\angle\lim_{\zeta\to 1}\frac{d(\tr_\v)_{\v(\zeta)} \cdot
F(\v(\zeta))}{\zeta-1}
\]
exists finitely, $A(\v,p)\in \R$ and $A(\v,p)\leq b$. Also, $b=\sup
A(\v,p)$, with the supremum taken as $\v$ varies among all Lempert's
projection devices $(\v, \tr_v)$ with $\v(1)=p$. \et

This result is contained in Theorem \ref{Julia continuous} and
Theorem \ref{bb}. One of the main ingredients in the proof is the
remarkable property that the projection of an infinitesimal
generator on every complex geodesic is still an infinitesimal
generator.

In order to give the proof of the previous results, in the first
section  we revise pluripotential theory in strongly convex domains
and in the second section we study iteration using the pluricomplex
Green function and the pluricomplex Poisson kernel. We should say
that, even if part of the results in section two are already known,
our present formulation seems to be new and, as we will prove later,
quite effective.  Section three is devoted to the interactions
between pluripotential theory and semigroups. In section four we
discuss a couple of examples on the boundary behavior of semigroups
and complete the proof of our characterization of BRFPs in terms of
the boundary behavior of the infinitesimal generator. As a
consequence, in Corollary \ref{boundary} we discuss {\sl stationary
points} of semigroups (namely those BRFPs for which the boundary
dilatation coefficient is less than or equal to $1$). In section
five we consider the non-linear resolvent of Reich and Shoikhet,
proving that every BRFP of the non-linear resolvent is a BRFP for
the semigroup (see Proposition \ref{non-linear}). Finally, in
section six we translate our results into the ball $\B^n$ where some
more explicit formulations, using automorphisms, are possible. In
this case, we also discuss the boundary behavior of the
infinitesimal generator at a BRFP under some boundness conditions
(see Corollary \ref{ball2}).

\medskip

Part of this work was done in Seville where the first named author
spent the entire month of March 2006. He wants to sincerely thank
the people at Departamento de Matem\'{a}tica Aplicada II at Escuela
Superior de Ingenieros in Universidad de Sevilla for the gentle
atmosphere and friendship he experienced there.

\section{Preliminary results on  pluripotential
theory in strongly convex domains}

For the definition, properties and further results about strongly
convex domains, we refer the reader to the nice monograph by Abate
\cite[Part 2]{Abate}. Likewise, for an introduction to
pluripotential theory with a special emphasis on  complex
Monge-Amp\`ere operators, we recommend the beautiful book by Klimek
\cite{Kl2} (a short introduction is also contained in
\cite{BrNotes}). Anyhow, for the sake of clearness, we are going to
give some basic definitions and define the tools we need later on.

\subsection{The pluricomplex Green function} Let $D\subset\subset \C^n$ be
a domain and $z\in D$. Define
\begin{equation*}
\mathcal K_{D,z}=\{u \ \hbox{plurisubharmonic in}\ D: u<0,
u(w)-\log\|z-w\|\leq O(1) \hbox{ as } w\to z\}.
\end{equation*}

The Klimek \cite{Kl} {\sl pluricomplex Green function} is defined as
\[
G_D(z,w):=\sup_{u\in \mathcal K_{D,z}}u(w).
\]
Such a function is plurisubharmonic in $D$, locally bounded in
$D\setminus\{z\}$ and has a logarithmic pole at $z$ (see \cite{Kl}
and \cite{Kl2}). If $D$ is hyperconvex (in particular, if $D$ is a
convex domain), then Demailly \cite{De} showed that $G_D$, extended
to be $0$ on $D\times \de D$, is continuous as a function
$G_D:D\times \overline{D}\to [-\infty,0)$. Moreover, from the work
of Lempert \cite{Le} and Demailly \cite{De}, it turns out that
$G_D(z,w)$ is the unique solution of the following homogeneous
Monge-Amp\`ere equation:
\[\begin{cases}
u\text{ plurisubharmonic in }D \\
(\partial \overline{\partial }u)^{n}=0 &\text{ in }D\setminus \{z\} \\
\lim_{w\rightarrow x}u(w)=0 &\text{ for all }x\in \partial D \\
u(w)-\log \left\vert w-z\right\vert =O(1) &\text{ as }w\rightarrow z.%
\end{cases}
\]

By the very definition, if $h:D\to D'$ is holomorphic, then for all
$z,w\in D$
\begin{equation}\label{pol}
G_{D'}(h(z),h(w))\leq G_D(z,w).
\end{equation}

In case $D$ is a bounded strongly convex domain with smooth boundary
(here and in the rest of the paper ``smooth'' means at least of
class $C^3$) Lempert \cite{Le} proved that $G_D(z,w)$ is smooth and
regular for $(z,w)\in \overline{D}\times \overline{D}\setminus {\sf
Diag}(\overline{D}\times \overline{D})$ and that
\begin{equation}\label{miol}
G_D(z,w)=\log\tanh k_D(z,w),
\end{equation}
where $k_D(z,w)$ is the Kobayashi distance of $D$ (for definition
and properties we refer to \cite{Abate} or to \cite{Kob}).

For instance, for $D=\D$ the unit disc in $\C$, the pluricomplex
Green function coincides with the usual (negative) Green function,
while for $D=\B^n$ the unit ball of $\C^n$ we have
\begin{equation}\label{Greenball}
G_{\B^n}(z,w)=\log \|T_z(w)\|,
\end{equation}
where $T_z:\B^n\to\B^n$ is any automorphism of $\B^n$ with the
property that $T_z(z)=0$.

\subsection{The pluricomplex Poisson kernel} Let $D\subset\subset \C^n$ be a strongly
convex domain with smooth boundary, $z_0\in D$ and let $p\in \de D$.
In the paper \cite{B-P}, Patrizio and the first quoted author
introduced a plurisubharmonic function $u_{D,p}:D\to (-\infty, 0)$
which extends smoothly on $\overline{D}\setminus\{p\}$ such that
$d(u_{D,p})_z\neq 0$ for all $z\in D$,  $u_{D,p}(q)=0$ for all
$q\in\de D\setminus\{p\}$ and $u_{D,p}$ has a simple pole at $p$
along non-tangential directions. Up to a real positive multiple, we
assume here that $u_{D,p}(z_0)=-1$. The function $u_{D,p}$ solves
the following homogeneous Monge-Amp\`ere equation:
\[\begin{cases}
u\text{ plurisubharmonic in }D \\
(\partial \overline{\partial }u)^{n}=0 &\text{ in }D \\
u<0 &\text{ in }D\\
u(w)=0 &\text{ for all }w\in \partial D\setminus \{p\} \\
u(w)\approx \left\Vert w-p\right\Vert ^{-1} &\text{ as }w\rightarrow
p \text{ non-tangentially}.
\end{cases}
\]

In the papers \cite{B-P} and \cite{BPT}, the authors prove that
$u_{D,p}$ shares many properties with the classical Poisson kernel
for the unit disk. In case $D=\D$ the unit disc in $\C$, the
function $u_{\D,p}$ (normalized so that $u_{\D,p}(0)=-1$) is in fact
the classical (negative) Poisson kernel. In case $D=\B^n$, the
pluricomplex Poisson kernel (normalized so that
$u_{\mathbb{B}^n,p}(0)=-1$) is given by
\[
u_{\mathbb{B}^n,p}(z)=-\frac{1-\|z\|^2}{|\la p-z, p\ra |^2}.
\]

The level sets of $u_{D,p}$ are exactly boundaries of Abate's
horospheres. Recall that a {\sl horosphere} $E_D(p,R)$ of center
$p\in \de D$ and radius $R>0$ (with respect to $z_0$) is given by
\[
E_D(p,R)=\{z\in D: \lim_{w\to p}[k_D(z,w)-k_D(z_0,w)]<
\frac{1}{2}\log R\}.
\]
Notice that the existence of the limit in the definition of
$E_D(p,R)$ is a characteristic of smooth strongly convex domains and
follows again from Lempert's theory (see
\cite[Theorem~2.6.47]{Abate}). Thanks to our normalization
$u_{D,p}(z_0)=-1$, it follows that
\begin{equation}\label{pluriP}
E_D(p,R)=\{z\in D: u_{D,p}(z)<-1/R\}.
\end{equation}

For the unit disk, these level sets are boundaries of horocycles
and, in case $D=\B^n$, these are boundaries of horospheres in $\B^n$
with center $p$, whose explicit expression is
\[
E_{\B^n}(p,R)=\{z\in \B^n: \frac{|1-\la z,p\ra|^2}{1-\|z\|^2}<R
\}.
\]

 More information about the properties of $u_{D,p}$ (such as
smooth dependence on $p$, extremality, uniqueness, relations with
the pluricomplex Green function, usage in representation formulas
for pluriharmonic functions) can be found in~\cite{BPT}.

\subsection{Lempert's projection devices}
We recall that a complex geodesic $\varphi: \D \to D$ is a
holomorphic isometry between $k_\D$ (the hyperbolic distance in
$\D$) and $k_D$. By Lempert's work (see \cite{Le} and
\cite{Abate}) given two points $z_0\in D$ and $z \in
\overline{D}$, there exists a unique complex geodesic $\varphi: \D
\to D$ such that $\varphi$ extends smoothly past the boundary,
$\varphi(0)=z_0$ and $\varphi(t)=z$, with $t\in (0,1)$ if $z \in
D$ and $t =1$ if $z \in
\partial D$. Moreover, for any such a complex geodesic there exists a
holomorphic retraction $\rho_\v : D \to \varphi(\D)$, {\em i.e.}
there exists a holomorphic map $\rho_\v:D\to D$  such that $\rho_\v
\circ \rho_\v =\rho_\v$ and $\rho_\v(z)=z$ for any $z
\in\varphi(\D)$.

Given a complex geodesic, there might exist many holomorphic
retractions to such geodesic, but the one constructed by Lempert
turns out to be the only one with affine fibers (see \cite[Section
3]{BPT}). We call such a $\rho_\v$ the {\em Lempert projection}
associated to $\varphi$.

Furthermore, we let $\tr_\v:= \varphi^{-1} \circ \rho_\v:D\to \D$
and call it the {\em left inverse } of $\varphi$, for $\tr_\v
\circ \varphi = {\sf id}_{\D}$. The triple $(\varphi, \rho_\v,
\tr_\v)$ is the so-called {\em Lempert projection device}.

For $D=\B^n$ the unit ball of $\C^n$ the image of the complex
geodesic through the points $z\neq w\in \overline{\B^n}$ is just
the {\sl one dimensional slice} $S_{z,w}:=\B^n\cap \{z+\zeta
(z-w): \zeta\in\C\}$. The Lempert  projection is thus given by the
orthogonal projection of $\B^n$ onto $S_{z,w}$.

By Lempert's very definition, if $\v:\D\to D$ is a complex geodesic,
then
\begin{equation}\label{GreenD-disco}
G_D(\v(\zeta),\v(\eta))=G_\D(\zeta,\eta)
\end{equation}
for all $\zeta,\eta\in\D$.

Finally, we mention \cite[p. 516]{B-P} that for any given Lempert
projection device $(\varphi, \rho_\v, \tr_\v)$ in $D$ with $\v(1)=p$
there exists $a_\v>0$ such that for all $\zeta\in\D$
\begin{equation}\label{PoissonD-disco}
u_{D,p}(\v(\zeta))=a_\v u_{\D,1}(\zeta).
\end{equation}

\section{Iteration theory by means of pluripotential theory}

Both the pluricomplex Green function and the pluricomplex Poisson
kernel can be used to describe dynamical properties of holomorphic
self-maps of a bounded strongly convex domain with smooth
boundary. The aim of this section is exactly to formulate the
results we need later on in terms of pluripotential theory.

All the results presented in this section are strongly based on
some known results about iteration  (mainly due to Abate, see
\cite{Abate}). However, for the aim of completeness, we sometimes
provide a sketch of some new direct proofs.

As a matter of notation, for a map $h:D\to D$ we denote by ${\sf
Fix}(h)$ the set of its fixed points in $D$, namely
\[
{\sf Fix}(h):=\{z\in D: h(z)=z\}.
\]

To begin with, we can reformulate a Schwarz-type lemma for strongly
convex domains as follows:

\bt\label{Schwarz} Let $D\subset\subset \C^n$ be a strongly convex
domain with smooth boundary and let $z_0\in D$. Let $h:D\to D$ be
holomorphic. Then $h(z_0)=z_0$ if and only if for all $z\in D$
\begin{equation}\label{iteraG}
G_D(z_0,h(z))\leq G_D(z_0, z).
\end{equation}
Moreover, if   equality holds in \eqref{iteraG} for some $z\neq z_0$
and $\v:\D\to D$ is the complex geodesic such that $\v(0)=z_0$ and
$\v(t)=z$ for some $t\in (0,1)$, it follows that $h\circ \v:\D\to D$
is a complex geodesic and $h:\v(\D)\to h(\v(\D))$ is an
automorphism.
 \et

\bpf The necessity and sufficiency of \eqref{iteraG} follows
directly from the very definition of $G_D$ and \eqref{pol}.

In order to prove the last statement, assume that $G_D(z_0,h(z))=
G_D(z_0, z)$ for some $z\in D$, $z\neq z_0$. Let $(\varphi,
\rho_\v, \tr_\v)$ be the Lempert projective device such that
$\v(0)=z_0$ and $\v(t)=z$ for some $t\in (0,1)$ and let $(\psi,
\rho_\psi, \tr_\psi)$ be the Lempert projective device such that
$\psi(0)=z_0$ and $\psi(r)=h(z)$ for some $r\in (0,1)$. Let
$\thh(\zeta):=\tr_\psi(h(\v(\zeta)))$ for $\zeta\in\D$. Notice
that $\thh(0)=0$. Then, $H(\zeta):=G_\D(0,\thh(\zeta))-G_\D(0,
\zeta)\leq 0$. Moreover, the function $\D\ni\zeta\mapsto H(\zeta)$
is harmonic on $\D\setminus\{0, \thh^{-1}(0)\}$ and bounded from
above, thus can be extended in a subharmonic way to all of $\D$.
We still call $H$ such an extension. For all $\zeta \in \D$, and
by \eqref{GreenD-disco}, it follows that
\[
G_\D(0,\thh(t))=G_D(\psi(0),\psi(\thh(t)))=G_D(z_0,h(z))=
G_D(z_0,z) =G_\D(0, t).
\]
By the maximum principle then $H(\zeta)\equiv 0$ and
$\thh(\zeta)=e^{i\theta}\zeta$ for some $\theta\in\R$, proving the
statement. \epf

\bd Let $D\subset\subset \C^n$ be a strongly convex domain with
smooth boundary, $p\in\de D$, and $h:D\to D$ holomorphic. The {\sl
boundary dilatation coefficient} $\al_h(p)\in (0,+\infty]$ is
defined as
\[
\al_h(p)=\inf_{q\in \de D}\{\sup_{z\in
D}\frac{u_{D,p}(z)}{u_{D,q}(h(z))}\}.
\]
\ed

As we show, this number can be characterized in several ways. Some
of them are widely used in the literature (see \cite{Abate} and
\cite{Br}). Indeed we have:

\bp\label{horosfera} Let $D\subset\subset \C^n$ be a strongly convex
domain with smooth boundary, $p\in \de D$ and $h:D\to D$
holomorphic. Then, the following are equivalent:
\begin{enumerate}
\item The boundary dilatation coefficient $\al_h(p)<+\infty$.

\item There exist a (necessarily unique) point $q\in \partial D$ and a number $\lambda>0$ such
that
\begin{equation}\label{med1}
h(E_D(p,R))\subseteq E_D(q,\lambda R), \text{ for all } R>0.
\end{equation}

\item It holds
\begin{equation}
\label{Marco} \frac{1}{2}\log \beta_h(p):=\liminf_{z\to
p}[k_D(z,z_0)-k_D(h(z),z_0)]<+\infty.
\end{equation}

\end{enumerate}
Moreover, if one of the statements holds, then
\begin{equation}
\beta_h(p)=\al_h(p)=\inf\{\lambda>0: \lambda \text{ satisfies}
\eqref{med1}\}.
\end{equation}
\ep

\begin{proof}
By the very definition (1) is equivalent to the existence of $q\in
\de D$ such that $u_{D,q}(h(z))\leq \frac{1}{\al_h(p)} u_{D,p}(z)$
for all $z\in D$. By \eqref{pluriP}, (1) and  (2) are equivalent
and  $\alpha_h(p)=\inf\{\lambda>0: \lambda \text{ satisfies}
\eqref{med1}\}$.

If (3) holds then  (2)  follows from Abate's version of the {\sl
Julia lemma} for strongly convex domains (see \cite[Theorem
2.4.16]{Abate}); also by the same token, $\beta_h(p)\geq
\inf\{\lambda>0: \lambda \text{ satisfies} \eqref{med1}\}$.

Finally, if (2) holds, let $\v:\D\to D$ be the complex geodesic such
that $\v(0)=z_0$ and $\v(1)=p$ and let $\tr_\v:D\to \D$ be its
left-inverse. Let $\th:=\tr_\v\circ h \circ \v:\D \to \D$. Since
$\v$ is an isometry between the Poincar\'e distance of $\D$ and the
Kobayashi distance of $D$ and $k_\D(\tr_\v(z),\tr_\v(w))\leq
k_D(z,w)$ for all $z,w\in D$, then it is easy to check that for all
$R>0$ it holds $\th(E_\D(1,R))\subseteq E_\D(1, \lambda R)$.
Therefore the classical Julia-Wolff-Carath\'eodory theorem implies
that $\beta_{\th}(1)<\infty$ and actually $\beta_{\th}(1)\leq
\lambda$. Now,
\begin{equation}\label{vaora}
\begin{split}
\frac{1}{2}\log\beta_h(p)&=\liminf_{w\to
p}[k_D(w,\v(0))-k_D(h(w),\v(0))]\\&\leq \liminf_{\zeta\to
1}[k_D(\v(\zeta),\v(0))-k_D(h(\v(\zeta)),\v(0))]\\
&\leq \liminf_{\zeta\to
1}[k_\D(\zeta,0)-k_\D(\tr_\v(h(\v(\zeta))),0)]=\frac{1}{2}\log
\beta_{\th}(1),
\end{split}
\end{equation}
which proves that $\beta_h(p)<+\infty$ and actually $\beta_h(p)\leq
\inf\{\lambda>0: \lambda \text{ satisfies } \eqref{med1}\}$, ending
the proof of the proposition.
\end{proof}

It is worth mentioning that by our very definition $\al_h(p)$ does
not depend on $z_0$, while {\sl a priori} the liminf in
\eqref{Marco} does. However, the independence of such liminf from
$z_0$ can be also shown directly, see   \cite[Lemma 6.1]{Br}.

We have the following version of Julia's lemma for strongly convex
domains:

\bt\label{Julia} Let $D\subset\subset \C^n$ be a strongly convex
domain with smooth boundary, $p\in \de D$ and $h:D\to D$
holomorphic. If the boundary dilatation coefficient
$\al_h(p)<+\infty$, then there exists a unique point $q\in \de D$
such that $h$ has non-tangential limit $q$ at $p$ and, for all
$z\in D$,
\begin{equation}\label{speedy}
u_{D,q}(h(z))\leq \frac{1}{\al_h(p)} u_{D,p}(z).
\end{equation}
Moreover, if  equality holds in \eqref{speedy} for some $z\in D$ and
$\v:\D\to D$ is the complex geodesic such that $\v(1)=p$ and
$\v(0)=z$, it follows that $h\circ \v:\D\to D$ is a complex geodesic
and $h:\v(\D)\to h(\v(\D))$ is an automorphism.
 \et

\bpf  By the very definition, if $\al_h(p)<+\infty$ then there
exists at least one $q\in \de D$ and a constant $C>0$ such that
\begin{equation}\label{speedy2}
u_{D,q}(h(z))\leq C u_{D,p}(z),
\end{equation}
for all $z\in D$. Since $u_{D,p}$ has a simple pole as $z\to p$
along non-tangential directions, the above inequality
\eqref{speedy2} implies that $h$ has non-tangential limit $q$ at
$p$. In particular, this implies that there exists {\sl at most
one} $q\in\de D$ such that $\sup_{z\in
D}\frac{u_{D,p}(z)}{u_{D,q}(h(z))}<+\infty$. Therefore,
\eqref{speedy} holds.

In order to prove the last statement, assume that $u_{D,q}(h(z))=
\frac{1}{\al_h(p)} u_{D,p}(z)$ for some $z\in D$. Let $(\varphi,
\rho_\v, \tr_\v)$ be the Lempert projective device such that
$\v(1)=p$ and $\v(0)=z$  and let $(\psi, \rho_\psi, \tr_\psi)$ be
the Lempert projective device such that $\psi(1)=q$ and
$\psi(0)=h(z)$. Write $\thh(\zeta):=\tr_\psi(h(\v(\zeta)))$ for
$\zeta\in\D$. By \eqref{PoissonD-disco} and \eqref{speedy} it
follows that $H(\zeta):=u_{\D,1}(\zeta)-\lambda
u_{\D,1}(\thh(\zeta))\leq 0$ for $\zeta\in \D$ and
$\lambda:=\al_h(p)a_\psi/a_\v$. The function $H$ is harmonic in
$\D$ and, by construction,
\[
H(0)=u_{\D,1}(0)-\lambda
u_{\D,1}(\thh(0))=\frac{1}{a_\v}u_{D,p}(z)-\frac{\lambda}{a_\psi}
u_{D,q}(h(z))=0.
\]
Thus the maximum principle implies that $H(\zeta)\equiv 0$, which
in turns implies that $\lambda=1$ and $\thh$ is the identity on
$\D$ and the statement follows. \epf

Let $D\subset \C^n$ be a bounded strongly convex domain with smooth
boundary, let $z_0\in D$ and let $p\in \de D$. Following Abate
(\cite{Abate}) we denote by $K(p,R)$ the {\sl K-region} with vertex
$p$ and radius $R>1$ defined as
\[
K(p,R)=\{z\in D: \lim_{w\to p}[k_D(z, w)-k_D(z_0,
w)]+k_D(z,z_0)<\log R\}.
\]
If $Q:D\to \C^n$ is a function, we write ${\sf K-}\lim_{z\to
p}Q(z)=L$ if for any sequence $\{z_k\}\subset D$ which tends to $p$
and belongs eventually to a K-region $K(p,R)$ for some $R>1$, it
follows $\lim_{k\to \infty} Q(z_k)=L$. Notice that if $Q$ has ${\sf
K-}$limit $L$ at $p$ then in particular it has non-tangential limit
$L$ at $p$.

\br\label{uso1} If $\al_h(p)<+\infty$ and $q\in\de D$ is the point
given by Theorem \ref{Julia} then actually $h$ has ${\sf K-}$limit
$q$ at $p$. This follows from Abate's version of the classical
Julia-Wolff-Carath\'eodory theorem, but also from \eqref{speedy},
since actually $u_{D,p}(z)\to -\infty$ when $z\to p$ inside a
K-region (see \cite[section 5]{BPT}).\er

The reason of the importance of  boundary dilatation coefficients in
iteration theory is that, while they give a global picture of the
dynamics of a self-map of $D$, they can be easily computed as radial
limits along any complex geodesic. We are going to state this fact
in a particular case which we need later. Before that we give the
following

\bd Let $D\subset \C^n$ be a strongly convex domain with smooth
boundary. Let $h:D\to D$ be holomorphic. We say that a point
$p\in\de D$ is a {\sl boundary regular fixed point}, BRFP for short,
if $h$ has non-tangential limit $p$ at $p$ and the boundary
dilatation coefficient $\al_h(p)<+\infty$. A BRFP with boundary
dilatation coefficient $\leq 1$ is also called a {\sl stationary
point}.  Likewise, those boundary regular fixed points with
$\al_h(p)>1$ are usually called {\sl boundary repelling fixed
points}.\ed

Now we can state the following version of Julia-Wolff-Carath\'eodory
theorem, due essentially to Abate:

\bt\label{cino} Let $D\subset\subset \C^n$ be a strongly convex
domain with smooth boundary. Let $h:D\to D$ be holomorphic and let
$p\in \de D$. Then $p$ is a BRFP for $h$ if and only if for
some---and hence any---Lempert projective device $(\varphi,
\rho_\v, \tr_\v)$ such that $\v(1)=p$ it follows
\begin{equation}\label{block}
\liminf_{(0,1)\ni r \to 1}\frac{|1-\tr_\v(h(\v(r)))|}{1-r}<+\infty.
\end{equation}
Moreover, if $p$ is a BRFP for $h$ then
\begin{equation}\label{block2}
\lim_{r\to 1
}\frac{1-\tr_\v(h(\gamma(r)))}{1-\tr_\v(\gamma(r))}=\al_h(p)
\end{equation}
for any curve $\gamma:[0,1)\to D$ such that
$\lim_{r\to1}\gamma(r)=p$, the curve in $\D$ given by $r\mapsto
\tr_\v(\gamma(r))$ converges non-tangentially to $1$ and $\lim_{r\to
1}k_D(\gamma(r), \rho_\v(\gamma(r)))=0$. In particular the map
$\tr_\v\circ h\circ \v:\D\to\D$ has BRFP at $1$ with boundary
dilatation coefficient $\al_h(p)$.
 \et

\bpf If $p$ is a BRFP for $h$ then the result follows from
\cite[Theorem 2.7.14]{Abate}.

Conversely, assume \eqref{block} holds. Then
\[
\liminf_{\zeta \to
1}\frac{1-|\tr_\v(h(\v(\zeta)))|}{1-|\zeta|}<\liminf_{(0,1)\ni r \to
1}\frac{|1-\tr_\v(h(\v(r)))|}{1-r}<+\infty.
\]
Thus the classical Julia-Wolff-Carath\'eodory theorem (see, {\sl
e.g.}, \cite{Abate}) implies that $1$ is a BRFP for $\zeta\mapsto
\tr_\v(h(\v(\zeta)))$ with  boundary dilatation coefficient
$a<+\infty$. Now, by \eqref{Marco}, taking into account that
$k_\D(\tr_\v(z),\tr_\v(w))\leq k_D(z,w)$ and arguing as in
\eqref{vaora} we find that
$\frac{1}{2}\log\al_h(p)\leq\frac{1}{2}\log a$, namely
$\al_h(p)<+\infty$. Theorem \ref{Julia} implies that $h$ has
non-tangential limit $q$ at $p$ for some $q\in \de D$. In order to
end the proof we need to show that $q=p$. To this aim, we first
notice that $\lim_{r\to 1} \tr_\v(h(\v(r)))=1$ forces $h(\v(r))$ to
tend  to $p$ as $r\to 1$ because
$\tr_\v(\overline{D}\setminus\{\v(\oD)\})\subset \D$ by
\cite[Proposition 1 p. 345]{Le2}. But $\v(\D)$ is transverse to $\de
D$ by Hopf's lemma and therefore $\v(r)\to p$ non-tangentially. This
implies that $\angle\lim_{z\to p}h(z)=p$ and we are done.
 \epf

In case a holomorphic self-map of $D$ has no fixed points in $D$,
 there always exists a particular stationary point (see
\cite[Theorem 2.4.23]{Abate}):

\bt[Abate]\label{wolfo} Let $D\subset\subset \C^n$ be a strongly
convex domain with smooth boundary. Let $h:D\to D$ be holomorphic.
If ${\sf Fix}(h)=\emptyset$ then there exists a unique point
$p\in\de D$, called the {\sl Denjoy-Wolff point} of $h$, such that
$p$ is a stationary point for $h$ and the sequence of iterates
$\{h^{\circ m}\}$ converges uniformly on compacta to the constant
map $D\ni z\mapsto p$. \et

Stationary points are quite special, as the following proposition
shows:

\bp\label{stationary} Let $D\subset\subset \C^n$ be a strongly
convex domain with smooth boundary. Let $h:D\to D$ be holomorphic.
Assume that $p\in \de D$ is a stationary point.
\begin{enumerate}
\item If ${\sf Fix}(h)\neq \emptyset$ then there exists a complex
geodesic $\v:\D\to D$ such that $\v(1)=p$ and $\v(D)\subseteq {\sf
Fix}(h)$. Moreover, for all $\theta\in\R$, the point
$\v(e^{i\theta})\in\de D$ is a stationary point for $h$ and
$\al_h(\v(e^{i\theta}))=1$. \item If ${\sf Fix}(h)= \emptyset$
then $p$ is the Denjoy-Wolff point of $h$ and $h$ has no other
stationary point in $\de D$.
\end{enumerate}
 \ep

\bpf  (1) Assume $z\in {\sf Fix}(h)$. Let $\v:\D\to D$ be the
complex geodesic such that $\v(0)=z$ and $\v(1)=p$. Consider the
holomorphic self-map of the unit disc $\psi(\zeta):=\tr_\v \circ h
\circ \v(\zeta)$. Then $\psi(0)=0$ and by Theorem \ref{cino}, $\psi$
has a stationary point at $1$. But then by the Herzig theorem
\cite{Hz} (see also the classical Wolff Lemma in \cite{Abate}) it
follows that $\psi(\zeta)\equiv \zeta$. Thus  for any $\zeta, \xi
\in \D$
\begin{equation*}\begin{split}
k_D(\varphi(\zeta), \varphi(\xi)) &\geq k_D(h(\varphi(\zeta)),
h(\varphi(\xi))) \geq k_D(\rho_\v(h(\varphi(\zeta))),
\rho_\v(h(\varphi(\xi))))\\&
  = k_D(\varphi(\psi (\zeta)), \varphi(\psi (\xi))) =
k_D(\varphi(\zeta), \varphi(\xi))=k_\D(\zeta,\xi),
\end{split}\end{equation*}
  forcing equality at all the
steps. In particular   $h \circ \varphi:\D\to D$ is a complex
geodesic such that $h(\v(0))=z$ and $h(\v(1))=p$. By the uniqueness
of complex geodesics passing through two given points of
$\overline{D}$ it follows that   $h \circ \varphi=\varphi$. Hence
$\v(\D)\subset {\sf Fix}(h)$.

Assertion (2) follows similarly. Indeed, let $q\in\de D$ be the
Denjoy-Wolff point of $h$. If $q\neq p$ then consider the complex
geodesic $\v:\D\to D$ such that $\v(-1)=q$ and   $\v(1)=p$ and  let
$\psi (\zeta):=\tr_\v \circ h \circ \v(\zeta)$. As before Theorem
\ref{cino} implies that $\psi$ has stationary points at $-1$ and
$+1$. Now the classical Wolff Lemma (see, {\sl e.g.}, \cite{Abate})
implies that $\psi (\zeta)\equiv \zeta$. Then we can proceed exactly
as before,  obtaining that $h(\v(\zeta))=\v(\zeta)$ for all
$\zeta\in \D$,  contradicting the hypothesis.   \epf

\section{Pluripotential theory and semigroups}

The aim of this section is to use the pluricomplex Green function
and the pluricomplex Poisson kernel to characterize infinitesimal
generators of semigroups of holomorphic self-maps of a strongly
convex domain and their dynamical properties.

We start recalling the following result (see \cite{A-S} for $D=\B^n$
and \cite[Theorem 2.5.24]{Abate}, \cite[Theorem A.1]{B-C-D} for the
general case)

\bt\label{fixwolf} Let $D\subset \C^n$ be a bounded strongly convex
domain with smooth boundary. Let $(\Phi_t)$ be a one-parameter
semigroup of holomorphic self-maps of $D$. Then
\begin{itemize}
\item either $\bigcap_{t\geq 0} {\sf Fix}(\Phi_t)\neq \emptyset$,
\item or ${\sf Fix}(\Phi_t)=\emptyset$ for all $t>0$, there exists a unique $\tau \in \partial D$ such that
$\tau$ is the Denjoy-Wolff point of $\Phi_t$ for all $t>0$ and there
exists $\beta\leq 0$ such that $\alpha_{\Phi_t}(\tau)=e^{\beta t}$.
\end{itemize}
\et

If a semigroup $(\Phi_t)$ has no fixed points in $D$, we call the
point $\tau\in\de D$ given by Theorem~\ref{fixwolf} the {\sl
Denjoy-Wolff point} of the semigroup.

\bd Let $D\subset \C^n$ be a bounded strongly convex domain with
smooth boundary. Let $(\Phi_t)$ be a   one-parameter semigroup of
holomorphic self-maps of $D$. A point $p\in \de D$ is called a
{\sl boundary regular fixed point for $(\Phi_t)$}, or a BRFP for
short, if $p$ is a BRFP for  $\Phi_t$ for all $t\geq 0$. The
family of boundary dilatation coefficients of $(\Phi_t)$ will be
denoted by $(\al_t(p))$. A BRFP for $(\Phi_t)$ for which
$\al_t(p)\leq 1$ for some $t>0$ is called a {\sl stationary point}
of the semigroup. \ed

The boundary dilatation coefficients at BRFP's form a semigroup in
$(\R_0^+,\cdot)$:

\bp\label{BRFPcomm} Let $D\subset \C^n$ be a bounded strongly
convex domain with smooth boundary. Let $(\Phi_t)$ be a
one-parameter semigroup of holomorphic self-maps of $D$. If
$p\in\de D$ is a BRFP for $(\Phi_t)$ then there exists $\beta\in
\R$ such that $\al_t(p)=e^{\beta t}$ for all $t\geq 0$. \ep

\bpf Let $(\varphi, \rho_\v, \tr_\v)$  be the Lempert projection
device associated to a complex geodesic such that $\varphi(1)=p$.
Consider the following family of functions $T_t:D\to \C$,
\[
T_t(z):=\frac{1-\tr_\v \circ \Phi_t (z)}{1-\tr_\v(z)}.
\]
By Theorem \ref{cino} it follows that $\lim_{(0,1)\ni r\to
1}T_t(\gamma(r))=\al_t(p)$ for any curve $\gamma:(0,1)\to D$ such
that $\lim_{r\to 1}\gamma(r)=p$, the curve $\tr_\v(\gamma(r))$
converges to $1$ non-tangentially and $k_D(\gamma(r),
\rho_\v(\gamma(r)))\to 0$ as $r\to 1$. By \cite[Proposition
3.4]{BrSNS}, it follows that $[0,1)\ni r \mapsto
\Phi_t(\varphi(r))$ satisfies the same three properties which are
satisfied by $\gamma$. Then for $s,t\geq 0$
 we have
\[
T_{t+s}(\varphi(r))=\frac{1-\tr_\v \circ \Phi_t
(\Phi_s(\varphi(r)))}{1-\tr_\v(\Phi_s(\varphi(r)))}\cdot\frac{1-\tr_\v
\circ
\Phi_s(\varphi(r))}{1-\tr_\v(\varphi(r))}=T_t(\Phi_s(\v(r))\cdot
T_s(\v(r)),
\]
and taking the limit as $r\to 1$ it follows that
$\alpha_{t+s}(p)=\alpha_t(p)\alpha_s(p)$. Since $\alpha_t(p)$ is
clearly measurable in $t$, this concludes the proof. \epf

Later we will see how the number $\beta$ in Proposition
\ref{BRFPcomm} can be computed using the infinitesimal generator of
the semigroup. Now we use the pluricomplex Green function to
characterize vector fields which are infinitesimal generators. For
this aim we need a lemma whose  simple proof is left to the reader:

\bl\label{convexidad-debil} Let $T>0$ be a positive real number and
let $g:[0,T]\rightarrow \R$ be a function such that
\begin{enumerate}
\item for all $a,b\in [ 0,T]$ and $\lambda \in [ 0,1]$ it holds
\begin{equation*}
g(\lambda a+(1-\lambda )b)\leq \max \{g(a),g(b)\};
\end{equation*}
\item there exists the (right-)derivative of $g$ at $0$ and $g'(0)>0.$
\end{enumerate}
Then $g$ is non-decreasing. \el

Now we can state and prove our characterizations of infinitesimal
generators:

\bt\label{caratterizza-Green} Let $D\subset\subset \C^n$ be a
strongly convex domain with smooth boundary. Let $F:D\to \C^n$ be
holomorphic. The following are equivalent:
\begin{enumerate}
\item The map $F$ is the infinitesimal generator of a semigroup
of holomorphic self-maps of $D$.
\item For all $z,w\in D$ with $z\neq w$ it follows that
\begin{equation}\label{decrescekd}
   d(k_D)|_{(z,w)}\cdot (F(z), F(w))\leq 0.
\end{equation}
\item For all $z,w\in D$ with $z\neq w$ it follows that
\begin{equation}\label{unoGb}
d(G_D)|_{(z,w)}\cdot (F(z), F(w))\leq 0.
\end{equation}
\item For each pair $z,w\in D,$ it follows
\begin{equation}
k_{D}(z-rF(z),w-rF(w))\geq k_{D}(z,w)
\end{equation}%
for all $r>0$ such that $z-rF(z)$ and $w-rF(w)$ belong to $D.$
\end{enumerate}
\et

 \bpf First of all we notice that by \eqref{miol} we have
$G_D(z,w)=\log \tanh k_D(z,w)$, and thus a simple computation shows
that (2) and (3) are equivalent.

Next, we claim that  (1) implies (3). Indeed, if $F$ is an
infinitesimal generator in $D$ and  $(\Phi_t)$ is the corresponding
semigroup generated by $F$ then by \eqref{pol}, for all $z,w\in D$
with $z\neq w$ it follows that for all $t\geq 0$
\[
G_D(\Phi_t(z), \Phi_t(w))- G_D(z,w)\leq 0
\]
and it is equal to zero for $t=0$. Computing the incremental ratio
in $t$ for $t=0$  we obtain~\eqref{unoGb}.

Now, assume (2) holds. For $w\in D$, consider the Cauchy problem
\[
\begin{cases}
\displaystyle{\frac{d \Phi}{d t}}=F \circ \Phi,\\
\Phi(0)=w
\end{cases}
\]
and denote by $\Phi_w:[0,\delta_w)\to D$ its maximal solution, for
some $\delta_w>0$. To show that $F$ is an infinitesimal generator,
it is enough to prove that for all $w$ it holds
$\delta_w=+\infty$.

To this aim, let $z,w\in D$ with $z\neq w$ and let
$\delta=\min\{\delta_z, \delta_w\}$. Let $g:[0,\delta)\ni t\mapsto
k_D(\Phi_z(t), \Phi_w(t))$. By uniqueness of solutions of the above
Cauchy problems, we know that, for all $t\in[0,\delta)$, we have
$\Phi_z(t)\neq \Phi_w(t)$.  According to Lempert's work \cite{Le},
\cite{Le2} (see also \cite[Proposition 2.6.40]{Abate}), the function
$g$ is smooth and differentiating with respect to $t$ we obtain by
\eqref{decrescekd}
\begin{equation*}
\begin{split}
g'(t)&=d(k_D)|_{(\Phi_z(t), \Phi_w(t))}\cdot (\frac{d\Phi_z(t)}{d
t}, \frac{d\Phi_w(t)}{d t})\\&=d(k_D)|_{(\Phi_z(t),
\Phi_w(t))}\cdot (F(\Phi_z(t)), F(\Phi_w(t)))\leq 0.
\end{split}
\end{equation*}
Therefore $g$ is non-increasing in $t$, namely
\begin{equation}\label{lop}
k_D(\Phi_z(t), \Phi_w(t))\leq k_D(\Phi_z(0), \Phi_w(0))=k_D(z,w).
\end{equation}
This implies that $\delta_z=\delta_w$ because, if for instance
$\delta_z<\delta_w$ then as $t\to \delta_z$ it would follow that
$\Phi_z(t)\to \de D$ while $\Phi_w(t)\to \Phi_w(\delta_z)\in D$,
and then $k_D(\Phi_z(t), \Phi_w(t))\to \infty$ contradicting
\eqref{lop}.

By the arbitrariness of $z, w\in D$, this means that {\sl for all}
$z\in D$ we have $\delta_z=\delta$. Hence, by well known results on
PDE's, we have a well defined analytic map $\Phi:D\times
[0,\delta)\to D$ which is holomorphic in $z\in D$ and such that
$\Phi(0,z)=z$ and $\frac{\de \Phi}{\de t}=F\circ \Phi$. Also,
$\Phi(t+s,z)=\Phi(t,\Phi(s,z))$ for all $s,t\geq 0$ such that
$s+t<\delta$ and $z\in D$. This implies that $\delta=\infty$.
Indeed, if $\delta<+\infty$, let $2\delta>t>\delta$ and let $s>0$ be
such that $t-s<\delta$, $s<\delta$. Define $\Phi_z(t):=\Phi(t-s,
\Phi(s,z))$. This is well defined and solve the Cauchy problem for
$z$, against the maximality of $\delta$.

Thus we have proved that (1), (2) and (3) are equivalent.

Now, let us prove that (4) implies (2). Let $z,w\in D,$ $z\neq w,$
and $r>0$ such that $z-rF(z)$ and $w-rF(w)$ belong to $D.$ By
convexity, $z-tF(z)$ and $w-tF(w)$ belong to $D$ for all $t\in [
0,r].$ Therefore, the function $g:[0,r]\rightarrow \R$ given by
\begin{equation*}
g(t)=k_{D}(z-tF(z),w-tF(w))
\end{equation*}
is well-defined and, again by Lempert's result, since $z\neq w,$ it
is differentiable at $0$. By hypothesis, $g(t)\geq g(0)$ for all
$t\geq 0$. Therefore $g'(0)\geq 0.$ But
\begin{equation*}
g'(0)=(dk_{D})|_{(z,w)}\cdot (-F(z),-F(w))=-(dk_{D})|_{(z,w)}\cdot
(F(z),F(w)).
\end{equation*}
Thus, $(dk_{D})|_{(z,w)}\cdot (F(z),F(w))\leq 0$, and (2) holds.

In order to finish the proof we show that (2) implies (4). To
proceed we consider  the following two possible cases:
\begin{itemize}
\item[I)] $(dk_{D})|_{(z,w)}\cdot (F(z),F(w))<0.$
\item[II)] $(dk_{D})|_{(z,w)}\cdot (F(z),F(w))=0$.
\end{itemize}

Case I). Fix $r>0$ such that $z-rF(z)$ and $w-rF(w)$ belong to $D.$
Then we have that $z-tF(z)$ and $w-tF(w)$ belong to $D$ for all
$t\in [0,r].$ Therefore, the function $g:[0,r]\rightarrow \R$ given
by
\begin{equation*}
g(t)=k_{D}(z-tF(z),w-tF(w))
\end{equation*}%
is well-defined and, since $z\neq w$, it is differentiable at $0$
with derivative given by
$g'(0)=(dk_{D})|_{(z,w)}\cdot (-F(z),-F(w))>0.$ Moreover, by \cite[Proposition 3.8]{R-S2},
given $%
z_{1},z_{2},w_{1},w_{2}\in D$ and $\lambda \in [ 0,1],$ we have that
\begin{equation*}
k_{D}(\lambda z_{1}+(1-\lambda )z_{2},\lambda w_{1}+(1-\lambda )w_{2})\leq
\max \{k_{D}(z_{1},w_{1}),k_{D}(z_{2},w_{2})\}.
\end{equation*}%
In particular, if $a,b\in [0,r]$ and $\lambda \in [0,1],$ we
have that
\begin{equation*}
\begin{split}
g(\lambda a&+(1-\lambda )b) =k_{D}(z-(\lambda a+(1-\lambda
)b)F(z),w-(\lambda a+(1-\lambda )b)F(w)) \\
&=k_{D}(\lambda (z-aF(z))+(1-\lambda )(z-bF(z)),\lambda
(w-aF(w))+(1-\lambda )(w-bF(w))) \\
&\leq \max \{k_{D}(z-aF(z),w-aF(w)),k_{D}(z-bF(z),w-bF(w))\} \\
&=\max \{g(a),g(b)\}.
\end{split}
\end{equation*}%
Therefore, $g$ satisfies the hypothesis of Lemma
\ref{convexidad-debil} and thus it is non-decreasing. Namely,
\begin{equation*}
k_{D}(z-tF(z),w-tF(w))\geq k_{D}(z,w)
\end{equation*}%
for all $t\in [0,r].$

Case II). Let $G:D\to \C^n$ holomorphic be an infinitesimal
generator in $D$ such that $(dk_{D})|_{(z,w)}\cdot (G(z),G(w))<0.$
Such a map can be constructed as follow. Up to translations we can
assume that $z=O$ the origin in $\C^n$. Let $a<0$. By convexity,
 the family of functions $\Phi_t:z\mapsto e^{at}z$ is a semigroup of
holomorphic self-maps of $D$. The associated infinitesimal generator
is $G(z)=az$. Therefore
\[
(dk_D)|_{(O,w)}\cdot (G(O), G(w))=a (dk_D)|_{(O,w)}\cdot (O, w).
\]
Now, the vector $(O,w)$ points outward with respect to the boundary
of the Kobayashi ball of center $O$ and radius $k_D(O,w)$ because
Kobayashi balls of convex domains are  convex (see, {\sl e.g.},
\cite[Proposition 2.3.46]{Abate}). Since Kobayashi balls are level
sets of $k_D$, this implies that $(dk_D)|_{(O,w)}\cdot (O, w)\neq
0$. Hence $d(k_D)|_{(O,w)}\cdot (G(O), G(w))\neq 0$ and, by the
already proved equivalence between (1) and (2), actually
$(dk_{D})|_{(z,w)}\cdot (G(z),G(w))<0.$

Now fix $\epsilon>0$ and consider the vector field $H:=F+\epsilon
G$. This is an infinitesimal generator of a semigroup of holomorphic
self-maps in $D$ (because $F+\epsilon G$ satisfies
\eqref{decrescekd} and by the equivalence between (1) and (2)). Now,
by construction, $(dk_D)|_{(z,w)}\cdot (H(z), H(w))<0$ and, for what
we proved in Case I), $k_{D}(z-rH(z),w-rH(w))\geq k_D(z,w)$ for all
$r>0$ such that $z-rH(z), w-rH(w)\in D$. Now, letting $\epsilon$
tends to $0$ we end the proof. \epf

As a corollary we have the following characterization of groups of
biholomorphisms of~$D$:

\bc  Let $D\subset\subset \C^n$ be a strongly convex domain with
smooth boundary. Let $F:D\to \C^n$ be holomorphic. The following are
equivalent:
\begin{enumerate}
\item The map $F$ is the infinitesimal generator of a group
of holomorphic self-maps of $D$.
\item For all $z,w\in D$ with $z\neq w$ it follows that
\begin{equation*}
   d(k_D)|_{(z,w)}\cdot (F(z), F(w))= 0.
\end{equation*}
\item For all $z,w\in D$ with $z\neq w$ it follows that
\begin{equation*}
d(G_D)|_{(z,w)}\cdot (F(z), F(w))= 0.
\end{equation*}
\end{enumerate}
 \ec
\bpf Apply Theorem \ref{caratterizza-Green} to $F$ and $-F$. \epf

\br\label{sciocco} In case $D=\B^n$ the unit ball of $\C^n$, using
\eqref{Greenball},  equation \eqref{unoGb} assumes a simple
expression given by
\begin{equation}\label{GGball}
\frac{\Re\la z, F(z)\ra}{1-\|z\|^2}+\frac{\Re\la w,
F(w)\ra}{1-\|w\|^2}\leq \Re \frac{\la F(z),w\ra +\la z, F(w)\ra
}{1-\la z, w\ra}.
\end{equation}
In fact, in case $D=\B^n$, Theorem \ref{caratterizza-Green} with
\eqref{GGball} replacing \eqref{unoGb}, was proven with different
methods by Reich and Shoikhet \cite[Theorem 2.1]{R-S}.
\er

For boundary regular fixed points, we have the following result:

\bt\label{Julia continuous} Let $D\subset\subset \C^n$ be a
strongly convex domain with smooth boundary. Let $F:D\to \C^n$ be
a holomorphic infinitesimal generator of a semigroup $(\Phi_t)$,
$\beta\in\R$ and $p\in\de D$. The following are equivalent:
\begin{enumerate}
\item The semigroup $(\Phi_t)$ has a BRFP at $p$ with boundary
dilatation coefficients $\al_t(p)\leq e^{\beta t}$ for all $t\geq
0$. \item $d( u_{D,p})_z\cdot F(z)+\beta u_{D,p}(z)\leq 0$ for all
$z\in D$.
\end{enumerate}
Moreover, if $p$ is a BRFP for $(\Phi_t)$ then the boundary
dilatation coefficient  of $\Phi_t$ is $\al_t(p)=e^{-tb}$ with
$b=\inf_{z\in D}d(u_{D,p})_z\cdot F(z)/u_{D,p}(z)$. \et

\bpf Suppose (1) holds. Then $u_{D,p}(\Phi_t(z))-e^{-t\beta}
u_{D,p}(z)\leq 0$ for all $t\geq 0$ and $z\in D$. In particular,
\begin{equation*}
\begin{split}
0\geq& \lim_{t\to 0^+}\frac{u_{D,p}(\Phi_t(z))-e^{-t\beta}
u_{D,p}(z)}{t}\\&=\frac{\de}{\de t}[u_{D,p}(\Phi_t(z))-e^{-t\beta}
u_{D,p}(z)]|_{t=0}=d(u_{D,p})_z\cdot F(z)+\beta u_{D,p}(z),
\end{split}
\end{equation*}
and (2) follows.

Conversely, assume (2) holds. Fix $z\in D$ and let
$g(t):=u_{D,p}(\Phi_t(z))-e^{-t\beta} u_{D,p}(z)$. We have to show
that $g(t)\leq 0$ for all $t\geq 0$. Deriving $g$, we obtain
\begin{equation*}
\begin{split}
g'(t)&= d(u_{D,p})_{\Phi_t(z)}\cdot\frac{\de\Phi_t}{\de
t}(z)+\beta e^{-\beta t}u_{D,p}(z)
\\&=d(u_{D,p})_{\Phi_t(z)}\cdot F(\Phi_t(z))+\beta e^{-\beta
t}u_{D,p}(z)\\&= d(u_{D,p})_{\Phi_t(z)}\cdot F(\Phi_t(z))+\beta
u_{D,p}(\Phi_t(z))-\beta g(t).
\end{split}
\end{equation*}
Therefore, using hypothesis (2) we have that for all $t\geq 0$
\begin{equation}\label{derivg}
g'(t)+\beta g(t) \leq 0.
\end{equation}
Now let $h(t):=-(g'(t)+\beta g(t))\geq 0$. Solving the
differential equation $g'(t)+\beta g(t)+h(t)=0$ with initial value
$g(0)=0$, we obtain
\[
g(t)=-e^{-\beta t}\int_0^t e^{\beta s} h(s) ds\leq 0,
\]
and thus (1) follows.

Finally, the last statement comes directly from Proposition
\ref{BRFPcomm} and the equivalence between (1) and (2).
  \epf

\br If $\beta\leq 0$ in Theorem \ref{Julia continuous}, it follows
that for all $z\in D$ the function $[0,+\infty)\ni t\mapsto
u_{D,p}(\Phi_t(z))-e^{-t\beta} u_{D,p}(z)$ is non-increasing.
Indeed, for $s>t$, we have
\begin{equation*}
\begin{split}
u_{D,p}(\Phi_s(z))-e^{-s\beta}
u_{D,p}(z)&=u_{D,p}(\Phi_{s-t}(\v_t(z)))-e^{-s\beta}
u_{D,p}(z)\\&\leq e^{-(s-t)\beta}[u_{D,p}(\Phi_t(z))-e^{-t\beta}
u_{D,p}(z)]\\&\leq u_{D,p}(\Phi_t(z))-e^{-t\beta} u_{D,p}(z).
\end{split}
\end{equation*}
 \er

We end up this section with a Berkson-Porta like characterization of
infinitesimal generators.

\bd Let $D\subset\subset \C^n$ be a strongly convex domain with
smooth boundary, $F:D\to \C^n$ holomorphic and $p\in \de D$. We
say that $F\in C_E^1(p)$ if for any horosphere $E_D(p,R)$ there
exists a $(n\times n)$-matrix $A$ such that
\[
\lim_{E_D(p,R)\ni z\to p} dF_z=A.
\]
\ed

\bt\label{charnobene} Let $D\subset\subset \C^n$ be a strongly
convex domain with smooth boundary. Let $p\in \de D$, $F:D\to
\C^n$ holomorphic and assume that $F\in C^1_E(p)$. Then $F$ is the
infinitesimal generator of a semigroup of holomorphic self-maps of
$D$ with a stationary point at $p$ if and only if
\begin{equation}\label{infh}
d(u_{D,p})_z\cdot F(z)\leq 0
\end{equation}
for all $z\in D$.
 \et

\bpf One direction follows directly from Theorem \ref{Julia
continuous}.

Conversely, assume that \eqref{infh} holds. Fix $w_0\in D$ and let
$\gamma:[0, \delta)\to D$ be the maximal solution of the Cauchy
problem
\begin{equation*}
\begin{cases}
\frac{d\gamma}{dt}=F\circ \gamma \\
\gamma(0)=w_0.
\end{cases}
\end{equation*}
It is enough to prove $\delta=+\infty$. Assume by contradiction
that $\delta<+\infty$. Let $g(t):=u_{D,p}(\gamma(t))$ for $t\in
[0,\delta)$. Deriving $g$, we obtain by \eqref{infh}
\[
g'(t)=d(u_{D,p})_{\gamma(t)}(\gamma'(t))=d(u_{D,p})_{\gamma(t)}(F(\gamma(t)))\leq
0.
\]
Thus for all $t\in [0,\delta)$ it follows $u_{D,p}(\gamma(t))\leq
u_{D,p}(\gamma(0))=u_{D,p}(w_0)$. This means that if $w_0\in
E_D(p,R)$ then $\gamma(t)$ belongs to  $E_D(p,R)$ for all
$t\in[0,\delta)$. In particular, since $\overline{E_D(p,R)}\cap
\de D=\{p\}$, it means that $\lim_{t\to \delta}\gamma(t)=p$.

Since $F\in C^1_E(p)$ and $\de E_D(p,R)$ is Lipschitz (it is
actually $C^{1,1}$ at $p$ and smooth elsewhere, see \cite[Section
4]{BPT}) by (a very simple form of) Whitney extension theorem there
exists a function $\tilde{F}:\C^n\to \C^n$ of class $C^1$ such that
$\tilde{F}|_{E_D(p,R)}=F$. If $\tilde{F}(p)=0$ then the Cauchy
problem
\begin{equation*}
\begin{cases}
\frac{d\eta}{dt}=\tilde{F}\circ \eta \\
\eta(\delta)=p
\end{cases}
\end{equation*}
has the unique solution $\eta(t)\equiv p$. In particular, $\gamma$
cannot reach $p$ in a finite time, which  gives us the searched
contradiction to $\delta<+\infty$.

To conclude the proof we are left to show that necessarily
$\tilde{F}(p)=0$. But this follows at once from the fact that for
any $z\in E_D(p,R)$ the solution of the Cauchy problem
\begin{equation*}
\begin{cases}
\frac{d\gamma^z}{dt}=\tilde{F}\circ \gamma^z \\
\gamma^z(0)=z
\end{cases}
\end{equation*}
is such that $\gamma^z(t)\in D$ for $t\in[0,\delta^z)$ for a
suitable $\delta^z\in (0,+\infty]$ and, arguing as for $\gamma$,
$\lim_{t\to \delta^z}\gamma^z(t)=p$.
 \epf

\br If $D=\D$ the unit disc in $\C$, then Theorem \ref{charnobene}
holds without any regularity assumption on $F$ at $p\in \de\D$.
Indeed, a direct computation shows that \eqref{infh} reduces
exactly to the Berkson-Porta formula \cite{Berkson-Porta}.

If $D=\B^n$ the unit ball of $\C^n$, a direct computation shows that
 \eqref{infh} corresponds to
\[
\frac{\Re \la F(z), z\ra}{1-\|z^2\|}\leq \Re \frac{\la F(z),
p\ra}{1-\la z,p\ra}.
\]
In fact, for $D=\B^n$ and with the additional hypothesis that
$F$ extends holomorphically through $\de\B^n$, Theorem
\ref{charnobene} follows from   \cite[Theorem 3.1]{AERS}.
 \er

\section{Boundary behavior of infinitesimal generators}

In all this section, $D$ denotes a bounded strongly convex domain
in $\C^n$ with smooth boundary.

Before proving the main result of this section, we examine two
significative examples. We will use a lemma whose proof can be
derived from the proof of \cite[Theorem 1.4]{B-C-D2}.

\begin{Lem}\label{lemmaCG} Let $a,b\in \C^n$ and $A\in \C^{n\times n}$
and
\[ G(z)=a-\la z,a\ra z-[Az+\la z,b\ra z].
\]
Then, $G$ is the infinitesimal generator of a continuous semigroup
of holomorphic self-maps of $\B^n$ if and only if
\begin{equation}\label{constrain}
|\la b,u\ra | \leq \Re \la Au, u\ra ,
\end{equation}
for all $u\in \partial \B^n$. Moreover, if equality holds at every
point of $\partial \B^n$, then $G$ is the infinitesimal generator
of a continuous group of holomorphic self-maps of $\B^n$.
\end{Lem}

\be\label{unos} Let us consider $F:\B^2\to \C^2$ given by
$F(z_1,z_2)=(0,-z_2/(1-z_1))$. Let $e_1=(1,0)$. By a direct
computation one can see that $F$ is an infinitesimal generator of a
semigroup $(\Phi_t)$ of holomorphic self-maps of $\B^2$ which
pointwise fixes the slice $\D\ni\zeta\mapsto (\zeta,0)$. Clearly
$d(u_{\B^2, e_1})_z \circ F(z)\leq 0$ for all $z\in \B^2$. Thus $F$
has a stationary point at $e_1$. Also, $\la F(z), e_1\ra =0$ for all
$z\in \B^2$ and therefore the radial limit $\lim_{(0,1)\ni r\to
1}\frac{F(re_1)}{1-r}=0$, as predicted by Theorem~\ref{Shoikhet}.
Now let us consider the following map
\[
\eta(z_1,z_2)=\left(\frac{-sz_2+(1-\beta)z_1+\beta}{-sz_2-\beta
z_1+1+\beta}, \frac{z_2+sz_1-s}{-sz_2-\beta z_1+1+\beta} \right)
\]
where $\Re \beta>0$ and $s=\sqrt{2\Re \beta}$. Notice that
$\eta:\B^2\to \B^2$ is a parabolic automorphism such that
$\eta(e_1)=e_1$ (see \cite[Example 5.1]{BisiBracci}). Hence
$\tilde{F}(z_1,z_2):=d\eta^{-1}_{\eta(z_1,z_2)}\cdot
F(\eta(z_1,z_2))$ is the infinitesimal generator of the semigroup
$(\tilde{\Phi}_t)$ of holomorphic self-maps of $\B^2$, where
$\tilde{\Phi}_t=\eta^{-1}\circ \Phi_t\circ \eta$. Thus
$\tilde{\Phi}_t$ pointwise fixes the slice $\eta^{-1}(\zeta,0)$
for $\zeta\in\D$. A direct computation shows that
$F(\eta(r,0))=(0,s)$ for $r\in(0,1)$. Thus, since
\[
d\eta_{re_1}=\frac{1}{(-r\beta +1+\beta)^2}\left(
              \begin{array}{cc}
                1 & s(r-1) \\
                s & \bar{\beta}(r-1)+1 \\
              \end{array}
            \right),
\]
it follows that $\tilde{F}(r,0)=(\beta(1-r)+1)(-s^2(r-1),s)$ for
$r\in (0,1)$. In particular the radial limit of $\tilde{F}$ at
$e_1$ is not zero. A direct computation shows that
\[
\lim_{(0,1)\ni r\to 1}\frac{\tilde{F}_1(re_1)}{r-1}=-s^2<0.
\]
Hence the semigroup of $\D$ generated by the ``projection''
$\tilde{F}_1(\zeta e_1)$ of $\tilde{F}$ to the slice $\D\ni
\zeta\mapsto \zeta e_1$ has Denjoy-Wolff point at $1$, with boundary
dilatation coefficients $(e^{-s^2 t})$. However, by construction,
$e_1$ is a stationary point for the semigroup $(\tilde{\Phi}_t)$,
with boundary dilatation coefficients all equal to $1$.
 \ee

\be\label{dues} Let $G:\B^2\to\C^2$ be the map defined by
$G(z_1,z_2)=(\frac{-1}{3}iz_1, iz_2)$. The map $G$ is an
infinitesimal generator of a (semi)group $(\Phi_t)$ of elliptic
automorphisms fixing the origin (see Lemma \ref{lemmaCG}). Let us
consider the automorphism of $\B^2$ given by
\[
\psi(z_1, z_2)=\left(\frac{-\sqrt{3} z_1}{2-z_2},
\frac{1-2z_2}{2-z_2} \right)
\]
(see \cite[Lemma 2.2.1]{Abate}). We note that $\psi=\psi^{-1}$.
Such an automorphism maps the slice $\{(z_1,z_2)\in\B^2: z_2=0\}$
to the slice $\{(z_1,z_2)\in\B^2: z_2=1/2\}$. Let us define
$P(z_1,z_2):=d\psi_{\psi^{-1}(z_1,z_2)}\cdot
G(\psi^{-1}(z_1,z_2))$. A direct computation shows that
\[
P(z_1,z_2)=\left(\frac{-2i}{3}z_1z_2,\frac{-i}{3}(z_2-2)(2z_2-1)
\right).
\]
Let $(\Phi^P_t)$ be the group of automorphisms generated by $P$.
Such a group is obtained by conjugation from the group $(\Phi_t)$
and therefore has only one fixed point. Since
$P(z_1,0)=(0,\frac{-2i}{3})$ (and $P$ is holomorphic past the
boundary) then $(\Phi^P_t)$ cannot have BRFP at any point of the
boundary of the slice $\{(z_1,z_2)\in\B^2: z_2=0\}$. However,
obviously
\[
\lim_{\D\ni \zeta\to 1}\frac{P_1(\zeta,0)}{\zeta-1}=0.
\]
Finally, let
\[
H(z_1,z_2)=\left( \frac{2i}{3}z_2(z_1-1),\frac{2i}{3}(1+z_2^2-z_1)
\right).
\]
The map $H:\B^2\to \C^2$ is the infinitesimal generator of a group
of automorphisms (see again Lemma \ref{lemmaCG}) with the property
that $H(e_1)=(0,0)$. Let us define $F=P+H$. Then $F$ is the
infinitesimal generator of a group of automorphisms and a direct
computation shows that
\[
F(z_1,z_2)=\left( -\frac{2i}{3}z_2,-\frac{5i}{3}z_2-\frac{2i}{3}z_1
\right).
\]
Thus the  group generated by $F$ has a unique fixed point at $(0,0)$
and no BRFPs on $\de \B^2$. However the semigroup generated by
$\D\ni\zeta\mapsto F_1(\zeta,0)$ on $\D$ is the trivial semigroup
and
\[
\lim_{\D\ni\zeta \to e^{i\theta}}\frac{F_1(\zeta,0)}{\zeta-1}=0
\]
for all $\theta\in \R$. Notice that the slice $\{(z_1,z_2)\in \B^2:
z_2=0\}$ contains the fixed point of the semigroup.
 \ee

The previous two examples show that, on the one hand, the
requirement that the radial limit exists in Theorem \ref{Shoikhet}
is sufficient but not necessary for the existence of BRFP's. Also,
even if a BRFP exists, say at $p\in\de\B^n$, the radial limit of
the incremental ratio of the projection of the infinitesimal
generator along $p$ might not give information on the boundary
dilatation coefficients of the semigroup at $p$. On the other
hand, the sole information on the existence of the limit of the
incremental ratio along a given point $p\in\de \B^n$ does not
imply existence of a BRFP at $p$. Last but not least, an
unexpected phenomenon takes place for infinitesimal generators:
the behavior of the semigroup generated by the ``restriction'' of
the infinitesimal generator to one complex geodesic---even a
complex geodesic containing fixed points of the semigroup---can be
completely different from the behavior of the semigroup in the
ball (cfr. Theorem \ref{cino}).

\bd Let $F:D\to \C^n$ be a holomorphic infinitesimal generator. For
a Lempert projection device $(\varphi, \rho_\v, \tr_\v)$ we will
denote by $f_\v(\zeta):=d(\tr_\v)_{\v(\zeta)}\cdot F(\v(\zeta))$ the
holomorphic vector field on $\D$. \ed

\bp\label{tecno-restr} Let $F:D \to \C^n$ be a holomorphic
infinitesimal generator. Let $(\v,\rho_\v, \tr_\v)$ be a Lempert
projection device with $p=\v(1)\in\de D$. Then the vector field
$f_\v(\zeta)$ is a holomorphic infinitesimal generator in $\D$.
Moreover, if there exists $\beta\in \R$ such that
$d(u_{D,p})_{z}\cdot F(z)+\beta u_{D,p}(z)\leq 0$ for all $z\in D$,
then $d(u_{\D,1})_{\zeta}\cdot f_\v(\zeta)+\beta u_{\D,1}(\zeta)\leq
0$ for all $\zeta\in \D$. \ep

\bpf Considering the pluricomplex Green function $G_D:D\times D \to
\R$. Its differential $dG_D:TD\times TD\to T\R$ can be decomposed as
$dG_D=d_zG_D+d_wG_D$ where, if $(u,v)\in TD\times TD$ we have
$dG(u,v)=d_zG_D(u)+d_wG_D(v)$. With this notation, Theorem
\ref{caratterizza-Green} implies  that for all $z\neq w$
\[
d(G_D)_{(z,w)}\cdot (F(z), F(w))=d_z (G_D)|_{(z,w)} \cdot F(z) +d_w
(G_D)|_{(z,w)} \cdot F(w)\leq 0.
\]
Now let $z=\v(\eta)$ and $w=\v(\zeta)$ for $\eta\neq \zeta\in \D$.
We claim that
\begin{equation}\label{restric}
d_w (G_D)|_{(\v(\eta),\v(\zeta))} \cdot F(\v(\zeta))=d_w
(G_D)|_{(\v(\eta),\v(\zeta))} (d(\rho_\v)_{\v(\zeta)}\cdot
F(\v(\zeta))).
\end{equation}
Assume that \eqref{restric} is true. According to
\eqref{GreenD-disco} we also have
$G_D(\v(\eta),\v(\zeta))=G_\D(\eta,\zeta)$ for all $\zeta\in\D$,
thus by \eqref{restric}
\begin{equation*}
\begin{split}
d_w (G_\D)|_{(\eta,\zeta)} (f_\v(\zeta))&=d_w
(G_D)|_{(\v(\eta),\v(\zeta))}\cdot (d\v_\zeta (f_\v(\zeta)))
\\&=d_w
(G_D)|_{(\v(\eta),\v(\zeta))}\cdot (d(\rho_\v)_{\v(\zeta)}
(F(\v(\zeta))))\\&=d_w (G_D)|_{(\v(\eta),\v(\zeta))} \cdot
F(\v(\zeta)).
\end{split}
\end{equation*}
A similar equation holds for $d_z
(G_D)|_{(z,w)}|_{(\v(\eta),\v(\zeta))} \cdot F(\v(\eta))$, swapping
the roles of $\eta$ and $\zeta$ in the previous argument. Thus
\[
d(G_\D)|_{(\eta,\zeta)}\cdot (f_\v(\eta),
f_\v(\zeta))=d(G_D)|_{(\v(\eta),\v(\zeta))}\cdot (F(\v(\eta)),
F(\v(\zeta))))\leq 0,
\]
for all $\zeta, \eta\in\D$ with $\zeta\neq \eta$, which implies that
$f_\v$ is an infinitesimal generator on $\D$ by
Theorem~\ref{caratterizza-Green}.

Now we are left to prove claim \eqref{restric}. Since $\rho_\v$ is
holomorphic, then $d \rho_\v= \de \rho_\v$ and the Lempert
projection $\rho_\v$ determines a holomorphic splitting of the exact
sequence of holomorphic bundles
\[
0\longrightarrow T\v(\D)\stackrel{\iota}{\longrightarrow}
TD|_{\v(\D)}\longrightarrow N_{\v(\D), D}\longrightarrow 0,
\]
 given by
\[
T_{\v(\zeta)}D=\iota(d\rho_\v(T_{\v(\zeta)}D))\oplus {\sf Ker}
d(\rho_\v)_{\v(\zeta)}.
\]
Now let $B_D(\v(0), R)=\{z\in D: k_D(z,\v(0))<R\}$ be a Kobayashi
ball for $D$ and let $\zeta_1\in \D$ be such that $\v(\zeta_1)\in
\de B_D(\v(0), R)$. It is known---and can be easily proven using
Lempert's special coordinates, see \cite{Le}, \cite{Le2}---that
\begin{equation}\label{tangent}
T^{\mathbb C}_{\v(\zeta_1)}\de B_D(\v(0), R)={\sf Ker}
d(\rho_\v)_{\v(\zeta_1)}
\end{equation}
where, as usual, $T^{\mathbb C}_{\v(\zeta_1)}\de B_D(\v(0), R)$
denotes the complex tangent space of $\de B_D(\v(0), R)$ at
$\v(\zeta_1)$. By \eqref{miol} it follows that
\[
\de B_D(\v(0), R)=\{z\in D: G_D(\v(0), z)=r\}
\]
for a suitable $r<0$, and in particular
\begin{equation}\label{kerG}
T_{\v(\zeta_1)}\de B_D(\v(0), R)={\sf Ker}(d_wG_D)_{(\v(0),
\v(\zeta_1))}.
\end{equation}

Now, consider the Kobayashi ball $B_D(\v(\eta), R)$ with
$R=R(\v(\zeta))>0$ such that $\v(\zeta)\in \de B_D(\v(\eta), R)$.
Equations \eqref{tangent} and \eqref{kerG} yield
\[
{\sf Ker} d(\rho_\v)_{\v(\zeta)}=T^{\mathbb C}_{\v(\zeta)} \de
B_D(\v(\eta), R)\subset T_{\v(\zeta)} \de B_D(\v(\eta), R)={\sf Ker}
(d_w G_D)_{(\v(\eta),\v(\zeta))}
\]
from which equation \eqref{restric} follows,  and the claim is
proved.

In order to prove the last assertion of the proposition, we argue
similarly as before. Let $\zeta\in \D$ and let $E_D(\v(1), R)$ be
the horosphere in $D$ which contains $\v(\zeta)$ on its boundary.
Then (again using Lempert's special coordinates, see \cite[p.
517]{B-P})
\begin{equation}\label{tangent2}
T^{\mathbb C}_{\v(\zeta_1)}\de E_D(\v(1), R)={\sf Ker}
d(\rho_\v)_{\v(\zeta_1)},
\end{equation}
and
\begin{equation}\label{kerP}
T_{\v(\zeta_1)}\de E_D(\v(1), R)={\sf
Ker}d(u_{D,\v(1)})_{\v(\zeta_1)}.
\end{equation}
Since $T^{\mathbb C}_{\v(\zeta)} \de E_D(\v(1), R)\subset
T_{\v(\zeta)} \de E_D(\v(1), R)$, equation \eqref{tangent2} yields
\begin{equation*}
d (u_{D,\v(1)})_{\v(\zeta)}\cdot F(\v(\zeta))=d
(u_{D,\v(1)})_{\v(\zeta)} (d(\rho_\v)_{\v(\zeta)}\cdot
F(\v(\zeta))),
\end{equation*}
and since $u_{D,\v(1)}\circ \v(\zeta)=a_\v u_{\D,1}(\zeta)$ for some
$a_\v>0$ by \eqref{PoissonD-disco}, we have
\begin{equation}\label{yap}
d(u_{\D,1})_{\zeta}\cdot f_\v(\zeta)+\beta
u_{\D,1}(\zeta)=\frac{1}{a_\v}[d(u_{D,p})_{\v(\zeta)}\cdot
F(\v(\zeta))+\beta u_{D,p}(\v(\zeta))].
\end{equation}
If $d(u_{D,p})_{z}\cdot F(z)+\beta u_{D,p}(z)\leq 0$ for all $z\in
D$  then $d(u_{\D,1})_{\zeta}\cdot f_\v(\zeta)+\beta
u_{\D,1}(\zeta)\leq 0$ for all $\zeta\in\D$ as stated. \epf

Before stating and proving the main result of this section we need a
preliminary lemma.

 \bl\label{MSP} Let $\D\subset \C$ be the unit disc
in $\C$. Let $G$ be the infinitesimal generator of a semigroup
$(\eta_t)$ of holomorphic self-maps of $\D$.  The following are
equivalent:
\begin{enumerate}
\item The point $1$ is a boundary regular fixed point for
$(\eta_t)$. \item There exists $C>0$ such that the radial limit
\[
\limsup_{(0,1)\ni r\to 1}\frac{|G(r)|}{1-r}\leq C.
\]
\end{enumerate}
Moreover, if $1$ is a BRFP for $(\eta_t)$ with boundary dilatation
coefficients $\al_t(1)=e^{b t}$ then
\[
\angle\lim_{\zeta\to 1}\frac{G(\zeta)}{\zeta-1}=b.
\]
\el \bpf If (1) holds then the result follows directly from
\cite[Theorem 1]{C-D-P}.

Conversely, hypothesis (2)  implies that $\lim_{r\to 1}G(r)=0$. By
Berkson-Porta's theorem~\ref{24}, there exists a point $b\in
\overline{\mathbb{D}}$ and a holomorphic function
$p:\mathbb{D\rightarrow C}$ with $\Re p\geq 0$ such that
\begin{equation*}
G(z)=(z-b)(\overline{b}z-1)p(z),\text{ \quad }z\in \mathbb{D}.
\end{equation*}
If $b=1$, we have that  $1$ is the Denjoy-Wolff point of the
semigroup $(\eta_t)$ and (1) follows. Otherwise, we have that
$\lim_{r\to 1}p(r)=0$. Then the function
$\varphi(z)=\frac{1-p(z)}{1+p(z)}$ is a self-map of the unit disc
and $\lim_{r\to 1}\varphi (r)=1$. By \cite[Proposition 4.13]{Po92},
the function $\varphi$ has angular derivative (possibly infinite) at
$1$. Thus, $p$, and so $G$, has angular derivative at $1$. That is,
there exists the radial limit
\[
\lim_{r\to 1}\frac{G(r)}{r-1}.
\]
By (2), such a limit is finite and again by \cite[Theorem
1]{C-D-P}, we obtain that 1 is a boundary regular fixed point of
the semigroup.\epf

\bt\label{bb} Let $D\subset\subset \C^n$ be a strongly convex
domain with smooth boundary, $F$ the infinitesimal generator of a
semigroup $(\Phi_t)$ of holomorphic self-maps of $D$, and $p\in
\de D$. The following are equivalent:
\begin{enumerate}
\item The point $p$ is a BRFP for $(\Phi_t)$. \item There exists
$C>0$ such that for any Lempert's projection device $(\v,\rho_\v,
\tr_\v)$ with $\v(1)=p$ it follows
\[
\limsup_{(0,1)\ni r\to 1}\frac{|f_\v(r)|}{1-r}\leq C.
\]
\end{enumerate}
Moreover, if $p$ is a BRFP for $(\Phi_t)$ with boundary dilatation
coefficients $\al_t(p)=e^{\beta t}$ then the non-tangential limit
\[
A(\v,p):=\angle\lim_{\zeta\to 1}\frac{f_\v(\zeta)}{\zeta-1}
\]
exists finite, $A(\v,p)\in \R$ and $A(\v,p)\leq \beta$. Also,
$\beta=\sup A(\v,p)$, with the supremum taken as $\v$ varies among
all  complex geodesics with $\v(1)=p$. \et

\bpf Suppose (1) holds. By Theorem \ref{Julia continuous} there
exists $\beta\in \R$ such that $d( u_{D,p})_z\cdot F(z)+\beta
u_{D,p}(z)\leq 0$ for all $z\in D$.  Let $(\v,\rho_\v, \tr_\v)$ be a
Lempert's projection device with $\v(1)=p$ and let $f_\v$ be the
associated vector field. By Proposition \ref{tecno-restr}, the map
$f_\v$ is an infinitesimal generator and satisfies $d(
u_{\D,1})_\zeta\cdot f_\v(\zeta)+\beta u_{\D,1}(\zeta)\leq 0$ for
all $\zeta\in\D$. By Theorem \ref{Julia continuous}, the semigroup
generated by $f_\v$ in $\D$ has a BRFP at $1$ with boundary
dilatation coefficients $\al_t(1)\leq e^{t\beta}$. By Lemma
\ref{MSP}, it follows that the non-tangential limit
$\angle\lim_{\zeta\to 1}\frac{f_\v(\zeta)}{\zeta-1}$ exists finite
and it is a real number less than or equal to $\beta$, and thus (1)
and part of the last statement are proved.

Suppose (2) holds. By Theorem \ref{Julia continuous}, it is enough
to show that there exists $\beta\in \R$ such that for all $z\in D$
it holds $d( u_{D,p})_z\cdot F(z)+\beta u_{D,p}(z)\leq 0$. Fix $z\in
D$ and let $\v:\D\to D$ be the complex geodesic such that $\v(0)=z$
and $\v(1)=p$. By Proposition \ref{tecno-restr}, the vector field
$f_\v$ is an infinitesimal generator in $\D$. Hypothesis (2) and
Lemma \ref{MSP} imply that $1$ is a BRFP for the semigroup generated
by $f_\v$ with boundary dilatation coefficients less than or equal
to $e^{Ct}$. Therefore Theorem~\ref{Julia continuous} applied to
$f_\v$ yields $d( u_{\D,1})_z\cdot f_\v(\zeta)+C u_{\D,1}(\zeta)\leq
0$ for all $\zeta\in\D$. By \eqref{yap} it follows that
$d(u_{D,p})_{\v(\zeta)}\cdot F(\v(\zeta))+C u_{D,p}(\v(\zeta))\leq
0$ for all $\zeta\in \D$ and thus, in particular for $\zeta=0$, we
have $d( u_{D,p})_{z}\cdot F(z)+C u_{D,p}(z)\leq 0$ as needed.

Finally, notice that, again by Theorem~\ref{Julia continuous}, the
previous arguments show also that if $p$ is a BRFP with boundary
dilatation coefficients $\al_t(p)=e^{\beta t}$ then $\beta$ is the
supremum of all $A(\v,p)$.
 \epf

\bc\label{boundary} Let $F:D\to \C^n$ be the holomorphic
infinitesimal generator of a semigroup $(\Phi_t)$ with a stationary
point $p\in\de\B^n$. Then for any Lempert's projection device
$(\v,\rho_\v, \tr_\v)$ with $\v(1)=p$
\begin{enumerate}
\item $\angle\lim_{\zeta\to 1} f_\v(\zeta) =0$, \item
$\angle\lim_{\zeta\to 1} f_\v(\zeta)/(\zeta-1)=A(\v,p)$ is a
finite real number, $A(\v,p)\leq 0$ and the boundary dilatation
coefficients $\al_t(p)=e^{t\beta}$ are such that $
A(\v,p)\leq\beta\leq 0$ for all $\v$.
\end{enumerate}
Moreover,
\begin{itemize} \item[a)] if $F(z)=0$ for some $z\in D$
then there exists a complex geodesic $\v:\D\to D$ with $\v(1)=p$
such that $F(\v(\zeta))=0$ for all $\zeta\in \D$ and all points of
$\v(\de \D)$ are stationary points for $(\Phi_t)$ with boundary
dilatation coefficients $\al_t(\v(e^{i\theta}))=1$  for all
$\theta\in\R$. Also  $A(\v,\v(e^{i\theta}))=1$  for all
$\theta\in\R$.
\item[b)] If $F(z)\neq 0$ for all $z\in D$ then $p$ is the
Denjoy-Wolff point of $(\Phi_t)$.
\end{itemize}
 \ec

\bpf  Taking into account that if $z\in D$ then $F(z)=0$ if and only
if $z\in {\sf Fix}(\Phi_t)$, the statement is a direct consequence
of Theorem \ref{bb} and Proposition \ref{stationary}. \epf

\section{Boundary repelling fixed points and  the non-linear resolvent}

In \cite{R-S} Reich and Shoikhet proved the following result.

\bt Let $D\subset\C^n$ be a bounded convex domain (not necessarily
strongly convex). Let $F:D\to \C^n$ be a holomorphic infinitesimal
generator of a semigroup $(\Phi_t)$ of holomorphic self-maps of
$D$. Then there exists a family $\{G_t\}$ of holomorphic self-maps
of $D$, with $G_0={\sf id}_D$, depending  on the parameter $t\in
[0,+\infty)$ such that for all $z\in D$ and $t\in [0,+\infty)$
\begin{equation}\label{non-resolvent}
G_t(z)-z=tF(G_t(z)),
\end{equation}
and for all $z\in D$
\begin{equation}\label{non-resolvent2}
F(z)=\lim_{t\to 0^+}\frac{G_t(z)-z}{t}.
\end{equation}
Moreover, if $z,w\in D$ are such that $w-z=tF(w)$ then $w=G_t(z)$.
\et

The above family $\{G_t\}$ is called the {\sl non-linear resolvent
of $F$}. This non-linear resolvent reads some dynamical properties
of the semigroup. Indeed, by \eqref{non-resolvent} (and
uniqueness) it follows easily that
\[
{\sf Fix}(G_t)=\{z\in D: F(z)=0\}={\sf Fix}(\Phi_t).
\]
In \cite[Proof of Corollary 1.6]{R-S}, it is also proved that if
$D\subset\subset\C^n$ is a strongly convex domain with smooth
boundary and $(\Phi_t)$ has no fixed points in $D$, and
$\tau\in\partial D$ is the Denjoy-Wolff point of the semigroup,
then
\[
G_t(E_D(\tau,R))\subseteq E_D(\tau,R),
\]
for all $R>0$ and all $t>0$. Then, by Theorem \ref{horosfera}, we
obtain that $\tau$ is a stationary point for all $G_t$. Since $F$
has no zeros in $D$, then ${\sf Fix} (G_t)=\emptyset$ and, by
Proposition \ref{stationary}, we conclude that $\tau$ is the
Denjoy-Wolff point of $G_t$, for all $t>0$. That is, the functions
$\Phi_t$ and $G_t$ share the same Denjoy-Wolff point.

For boundary regular fixed points, we can prove:

\bp\label{non-linear} Let $D\subset\subset \C^n$ be a strongly
convex domain with smooth boundary. Let $F:D\to \C^n$ be a
holomorphic infinitesimal generator with associated semigroup
$(\Phi_t)$, non-linear resolvent $\{G_t\}$, and $p\in\de D$.
Suppose there exists $\beta\in \R$ such that for any $t>0$ the
point $p$ is a BRFP for $G_t$ with boundary dilatation
coefficients $\al_{G_t}(p)\leq e^{t\beta}$. Then $p$ is a BRFP for
$(\Phi_t)$ with boundary dilatation coefficients $\al_t(p)\leq
e^{t\beta}$. \ep

\bpf By Theorem \ref{Julia continuous} it is enough to prove that
$d( u_{D,p})_z\cdot F(z)+\beta u_{D,p}(z)\leq 0$ for all $z\in D$.
Fix $z\in D$. Since $p$ is a BRFP for $\{G_t\}$ and
$\al_{G_t}(p)\leq e^{t\beta}$, by Theorem \ref{Julia}, we have
\[
g(t):=u_{D,p}(G_t(z))-e^{-t\beta}u_{D,p}(z)\leq 0
\]
for all $t\in [0,\infty)$.  Since $g(0)=0$, using
\eqref{non-resolvent} and the fact that $G_t(z)\to z$ for $t\to
0^+$ by \eqref{non-resolvent2}, we have
\begin{equation*}\begin{split}
0&\geq \lim_{t\to 0^+}\frac{g(t)}{t}=\lim_{t\to
0^+}\frac{u_{D,p}(G_t(z))-u_{D,p}(z)}{t}+\beta
u_{D,p}(z)\\&=\lim_{t\to
0^+}\frac{u_{D,p}(tF(G_t(z))+z)-u_{D,p}(z)}{t}+\beta
u_{D,p}(z)\\&=d(u_{D,p})_z\cdot F(z) +\beta u_{D,p}(z),
\end{split}\end{equation*}
proving the statement.
 \epf

\br If $p\in \de D$ is a BRFP for $\{G_t\}$ with boundary
dilatation coefficients $\al_{G_t}(p)\leq e^{\beta t}$ for some
$\beta\in\R$, then (see Remark \ref{uso1}), for any $t\in
[0,\infty)$, it follows that ${\sf K-}\lim_{z\to p}G_t(z)=p$. In
particular, from \eqref{non-resolvent}, it follows that for all
$t>0$ and $R>1$
\[
\lim_{G_t(K(p,R))\ni z\to p} F(z)=0.
\]
However, even in this case, $F$ might not have radial limit $0$ at
$p$. In fact, looking at the infinitesimal generator $\tilde{F}$ in
Example \ref{unos}, one  easily sees that the non-linear resolvent
$\{G_t\}$ has a BRFP at $e_1$ with boundary dilatation coefficients
$\al_{G_t}(e_1)=1$ (because by construction $G_t(z)=z$ on a complex
geodesic containing $e_1$ on its boundary). But $\tilde{F}$ does not
have radial limit $0$ at $e_1$.  \er

The converse to Proposition \ref{non-linear} is false, as the
following example shows:

\be Let $f(\zeta)=1-\zeta^2$. Then $f$ is the infinitesimal
generator of a group of hyperbolic automorphisms in $\D$, with
Denjoy-Wolff point $1$ and boundary repelling fixed point $-1$. It
is easy to check that the non-linear resolvent of $f$ is given by
$$
G_t(z)=
\frac{1}{2t}(-1+\text{exp}(\frac{1}{2}\text{Log}(1+4t(t+z)))
$$
for all $t>0$ and $z\in \D$. Now, a direct computation shows
$G_t(-1)=\frac{|2t-1|-1}{2t}\neq -1$. \ee

\section{Boundary behavior in the the unit ball}

In this section we translate our  results on BRFP's for semigroups
of the unit ball $\B^n\subset\C^n$, where most expressions have
computable  forms.

In order to simplify our statements and without loss of
generality, we will assume that, up to conjugation, the base point
is $e_1=(1,0,\ldots,0)\in \de \B^n$.

\bt\label{ball1} Let $F:\B^n\to \C^n$ be the infinitesimal
generator of a semigroup $(\Phi_t)$ of holomorphic self-maps of
$\B^n$. The following are equivalent:
\begin{enumerate}
\item The point $e_1\in\de\B^n$ is a BRFP for $(\Phi_t)$. \item
There exists $C>0$ such that for all automorphisms $H=(H_1,\ldots,
H_n):\B^n\to \B^n$ such that $H(e_1)=e_1$ it follows
\begin{equation}\label{condB}
\limsup_{(0,1)\ni r\to 1} \frac{|
d(H_1)_{H^{-1}(re_1)}(F(H^{-1}(re_1)))|}{1-r}\leq C.
\end{equation}
\end{enumerate}
Moreover, if $e_1$ is a BRFP for $(\Phi_t)$ with boundary dilatation
coefficients $\al_t(e_1)=e^{\beta t}$ then the non-tangential limit
\[
A(H,e_1):=\angle\lim_{\zeta\to 1}\frac{ d(H_1)_{H^{-1}(\zeta
e_1)}(F(H^{-1}(\zeta e_1))) }{\zeta-1}
\]
exists finitely, $A(H,e_1)\in \R$ and $A(H,e_1)\leq \beta$. Also,
$\beta=\sup A(H,e_1)$, with the supremum taken as $H$ varies among
all  automorphisms of $\B^n$ with $H(e_1)=e_1$.
 \et

\begin{proof}
The result follows from Theorem \ref{bb} as soon as one realizes
how Lempert's projection devices  in the unit ball are related to
automorphisms of $\B^n$. Indeed, thanks to the double transitivity
of the group of automorphisms of $\B^n$ on $\de \B^n$, any complex
geodesic $\v:\D\to\B^n$ of $\B^n$ passing through $e_1$ can be
written as $\zeta\mapsto H^{-1}(\zeta e_1)$ for some suitable
automorphism $H:\B^n\to \B^n$ . The associated Lempert projection
$\rho_\v$ is thus given by $\rho_\v(z)=H^{-1}(\la H(z),e_1\ra
e_1)=H^{-1}(H_1(z),0\ldots,0)$ and the left inverse is
$\tr_\v(z)=H_1(z)$. Therefore
\[
f_\v(\zeta)=d(\tr_\v)_{\v(\zeta)}\cdot
F(\v(\zeta))=d(H_1)_{H^{-1}(\zeta e_1)}\cdot F(H^{-1}(\zeta e_1)),
\]
from which the statement follows.
\end{proof}

In the statement of Theorem \ref{ball1}, the sufficient condition
for $e_1$ to be a BRFP can be checked  considering only the class of
{\sl parabolic} automorphisms $H$ (namely, those for which the
boundary dilatation coefficient at $e_1$ is $1$). For the sake of
clearness, we examine in detail  the case $n=2$. In such a case we
can limit ourselves to (parabolic) automorphisms of $\B^2$ of the
form
\begin{equation}\label{Hauto}
H_{s,\theta}(z)=\frac{(-sz_2+(1-\beta)z_1+\beta,
e^{i\theta}(z_2+sz_1-s))}{-sz_2-\beta z_1+1+\beta},
\end{equation}
where $  \beta\geq 0$, $s=\sqrt{2\beta}$ and $\theta\in\R$. Notice
that $(H_{s,\theta})^{-1}$ also has the same form of $H_{s,\theta}$
and the differential of $H_{s,\theta}$ at $e_1$ is
\[
d(H_{s,\theta})_{e_1}=\left(
           \begin{array}{cc}
             1 & 0 \\
             se^{i\theta} & e^{i\theta} \\
           \end{array}
         \right).
\]
If $\v:\D\to \B^2$ is a complex geodesic with $\v(1)=e_1$ and we
write $\v'(1)$ in projective coordinates as
$\v'(1)=[1:se^{i\theta}]$, with $s\geq 0$ and $\theta\in\R$, the
corresponding $H_{s,\theta}$ in \eqref{Hauto} is such that
$H_{s,\theta}(\D\times\{0\})=\v(\D)$ and therefore, by uniqueness of
complex geodesics, $\v(\zeta)=H_{s,\theta}(\psi(\zeta),0)$ for some
automorphism $\psi$ of $\D$. Thus, in the statement of Theorem
\ref{ball1} for $n=2$, it is enough to check condition \eqref{condB}
for $H$ belonging to the class of $H_{s,\theta}$'s.

Theorem \ref{ball1} and the previous observation   can be used to
obtain the boundary behavior of   infinitesimal generators with some
bounds on the image. To explain this fact, we prove the following
corollary in $\B^2$, which can be easily generalized to $\B^n$ for
any $n\geq 2$, and can be considered a Julia-Wolff-Carath\'eodory
type theorem for infinitesimal generators.

\bc\label{ball2} Let $F:\B^2\to \C^2$ be the infinitesimal generator
of a semigroup with a BRFP at $e_1$. Suppose there exist  a
horosphere $E_{\B^2}(e_1,R)$ and two distinct points $a_0,a_1\in \C$
such that $F_1(E_{\B^2}(e_1,R)) \subset \C\setminus\{a_0,a_1\}$.
Then
\begin{enumerate}
\item $F_1$ has non tangential limit $0$ at $e_1$, namely, $\angle
\lim_{z\to e_1}F_1(z)=0$.
\item $\angle \lim_{\zeta\to 1} (1-\zeta) F_2(\frac{((1-\beta)\zeta+\beta,
e^{i\theta}( s\zeta-s))}{-\beta \zeta+1+\beta})=0$ for all
$\beta\geq 0$, $s=\sqrt{2\beta}$ and $\theta\in \R$.
\end{enumerate}
\ec
\begin{proof}
By the very definition of horospheres in $\B^2$, there exists a ball
  $B\subset E_{\B^2}(e_1,R)$ such that $B$ is tangent to $\B^2$ at
$e_1$. Let $\{z_k\}\subset\B^2$ be any sequence converging to $e_1$
non-tangentially. Then the sequence $\{z_k\}$ is eventually
contained in $B$. Hence
\begin{equation}\label{cilp}
k_{E_{\B^2}(e_1,R)}(z_k, \la z_k, e_1\ra e_1)\leq k_B(z_k, \la z_k,
e_1\ra e_1).
\end{equation}
For $k\to \infty$ we have $k_B(z_k, \la z_k, e_1\ra e_1)\to 0$
because $z_k\to e_1$ non-tangentially in $B$ (and non-tangential
sequences are special in the sense of Abate \cite[Lemma
2.2.24]{Abate}). Therefore
\begin{equation}\label{cilp2}
\lim_{k\to \infty} k_{E_{\B^2}(e_1,R)}(z_k, \la z_k, e_1\ra e_1)=0.
\end{equation}
Now let $g :=F_1|_{E_{\B^2}(e_1,R)}:E_{\B^2}(e_1,R)\to \mathcal
L:=\C\setminus\{a_0, a_1\}$. By the monotonicity of Kobayashi
distance we have
\[
\omega_{\mathcal L}(g(z_k),g(\la z_k, e_1\ra e_1))\leq
k_{E_{\B^2}(e_1,R)}(z_k, \la z_k, e_1\ra e_1),
\]
and   \eqref{cilp2} forces
\[
\lim_{k\to \infty} \omega_{\mathcal L}(g(z_k),g(\la z_k, e_1\ra
e_1))=0.
\]
Since $\mathcal L$ is hyperbolic, this means that if $g(\la z_k,
e_1\ra e_1)$ tends to some $b\in \C$ then $g(z_k)$ must have the
same limit as $k\to \infty$. By \eqref{condB} it follows that
$g(\zeta e_1)$ has non-tangential limit $0$ at $1$. Since $z_k\to
e_1$ non-tangentially,  the same does   $\{\la z_k, e_1\ra e_1\}$.
Then
\[
\lim_{k\to \infty}F_1(z_k)=\lim_{k\to \infty}g(z_k)=\lim_{k\to
\infty}g(\la z_k, e_1\ra e_1)=0,
\]
proving that $F_1$ has non-tangential limit $0$ at $e_1$.

As for (2),  from Theorem \ref{ball1} with $H=H_{s',\theta'}$ (for
$s'\geq 0$, $\theta'\in\R$) as in \eqref{Hauto}, we have
\[
\angle\lim_{\zeta \to
1}d((H_{s',\theta'})_1)_{H_{s',\theta'}^{-1}(\zeta
e_1)}(F(H_{s',\theta'}^{-1}(\zeta e_1)))=0.
\]
By the very definition of $H_{s',\theta'}$ (and keeping in mind that
$(H_{s',\theta'})^{-1}=H_{s,\theta}$  for some $s\geq 0$ and
$\theta\in\R$) an easy computation shows that
\[
d((H_{s',\theta'})_1)_{H_{s',\theta'}^{-1}(\zeta
e_1)}(F(H_{s',\theta'}^{-1}(\zeta
e_1)))=C(\zeta)[F_1(H_{s,\theta}(\zeta e_1))+
(1-\zeta)F_2(H_{s,\theta}(\zeta e_1))],
\]
where $C(\zeta)$ is a smooth function which tends to some real
number $C\neq 0$ for $\zeta\to 1$. Thus, since $F_1$    has
non-tangential limit $0$ at $e_1$ by (1), statement (2) follows.
\end{proof}

\be The infinitesimal generator $F(z_1,z_2)=(0,-z_2/(1-z_1))$ in
Example \ref{unos} has the  boundary behavior prescribed by
Corollary \ref{ball2} at $e_1$. Notice that $F_2$ has not
(non-tangential) limit $0$ at $e_1$. \ee

\end{document}